\documentclass[final]{ectaart}
\usepackage[OT1]{fontenc}
\usepackage{graphicx}
\usepackage{amssymb}
\usepackage{amsthm}
\usepackage[cmex10]{amsmath}
\usepackage{natbib}
\usepackage{multirow}
\usepackage[colorlinks,citecolor=blue,urlcolor=blue,unicode,verbose]{hyperref}

\DeclareMathAlphabet\scr{U}{scr}{m}{n}
\DeclareMathAlphabet\scr{U}{scr}{m}{n}
\SetMathAlphabet\scr{bold}{U}{scr}{b}{n}
\DeclareFontFamily{U}{scr}{\skewchar\font'177}
\DeclareFontShape{U}{scr}{m}{n}{<-6>rsfs5<6-8>rsfs7<8->rsfs10}{}
\DeclareFontShape{U}{scr}{b}{n}{<-6>rsfs5<6-8>rsfs7<8->rsfs10}{}
\newtheorem{theorem}{Theorem}[section]

\newtheorem{corollary}[theorem]{Corollary}

\newtheorem{definition}[theorem]{Definition}

\newtheorem{lemma}[theorem]{Lemma}

\newtheorem{proposition}[theorem]{Proposition}
\newtheorem{remark}[theorem]{Remark}

\numberwithin{equation}{section}
\theoremstyle{plain}
\newtheorem{assumption}[theorem]{Assumption}
\newtheorem{convention}[theorem]{Convention}
\renewenvironment{proof}[1][\sc Proof]{\noindent\textbf{#1.} }{\hfill Q.E.D.}
\def\softl{l\kern-0.35ex\raise0.1ex\hbox{'}\kern-0.15ex}
\def\qed{\hfill\hbox{$\square$}}
\newcommand{\hookrightdoublearrow}{  \hookrightarrow\mathrel{\mspace{-15mu}}\rightarrow}
\newcommand{\conv}{\mathbin{\square}}

\newcommand{\cl}{\operatorname{cl}}
\begin{document}

\begin{frontmatter}
\title{Convex duality and Orlicz spaces in expected utility maximization\protect\thanksref[*]{T1}}
\runtitle{Convex duality for expected utility maximization}
\thankstext[*]{T1}{We would like to thank Jan Kallsen for his kind hospitality in September
2014 and Teemu Pennanen for guidance in all matters conjugate. We also thank Martin
Cripps, Christoph K\"uhn, and two anonymous referees for helpful comments. The paper has benefited significantly
from discussions with Christa Cuchiero, Irene Klein and Josef Teichmann on
FTAP and supermartingale deflators. Last but not least, we are grateful to
Fabio Maccheroni and Massimo Marinacci for putting us in touch to work on
this topic back in 2009.}
\begin{aug}
\author{\fnms{Sara} \snm{Biagini}\ead[label=e1]{sbiagini@luiss.it}}
\author{\fnms{Ale\v{s}} \snm{\v{C}ern\'{y}}\ead[label=e2]{ales.cerny.1@city.ac.uk}}

\runauthor{S. Biagini and A. \v{C}ern\'{y}}

\address{LUISS Guido Carli, Rome, \printead{e1}}

\address{Cass Business School, City, University of London, \printead{e2}}
\end{aug}

\begin{abstract}
In this paper we report further progress towards a complete theory of
state-independent expected utility maximization with semimartingale price
processes for arbitrary utility function. Without any technical assumptions
we establish a surprising Fenchel duality result on conjugate Orlicz spaces,
offering a new economic insight into the nature of primal optima and
providing fresh perspective on the classical papers of \citet{kramkov.schachermayer.99,kramkov.schachermayer.03}. 
The analysis points to an intriguing interplay
between no-arbitrage conditions and standard convex optimization and
motivates study of the Fundamental Theorem of Asset Pricing (FTAP) for
Orlicz tame strategies.
\end{abstract}


\end{frontmatter}

\section{Introduction}\label{sect: intro}

Utility maximization is a fundamental tenet of normative economic theory
and, as its most classical embodiment, \textquotedblleft expected utility
remains the primary model in numerous areas of economics dealing with risky
decisions\textquotedblright\ \citep{moscati.16}. Although a rigorous
axiomatic foundation of expected utility appeared early 
\citep{neumann.morgenstern.44} there remains a long-standing open problem in
the theoretical description of expected utility maximization in a purely
financial dynamic stochastic setting. Our aim is to offer new insights in
this direction.

The paper studies the mechanics of wealth transfer from initial date $0$ to
some terminal date $T$. A single agent whose preferences over terminal
wealth are represented by expected utility under given subjective
probability $P$ decides, continuously in time, how to allocate her wealth
among one risk-free and finitely many risky assets modeled by a
semimartingale price process $S$. There is no intermediate consumption and
no production or labour income. The main concern of the paper is finding a
suitable class of trading strategies that makes the problem well-defined. 
This is a non-trivial task because, as observed by \citet{harrison.kreps.79}, 
 unrestricted trading in continuous time permits so-called doubling 
strategies that create something out of nothing with certainty even 
when trading on a martingale. 

In this paper we make three distinct contributions to the literature. Firstly, the
Orlicz space framework unifies different strands of currently fragmented
literature on utility maximization and absence of arbitrage. Coupled with
convex duality it also conveys strong economic intuition. The unifying
framework, its economic interpretation, and links to the relevant literature are
presented in Sections \ref{sect: intro utility}-\ref{sect: intro Fenchel}.

Our second contribution is a new Fenchel duality result (Theorem~\ref{thm: strong duality}) 
which allows us to remove singular parts in the dual problem and offer a new
interpretation of the resulting duality as an `effective completion' of the
market. Effective completion means that the complete market represented by
the dual optimizer does not contain the entirety of the original opportunity
set but only those elements that have finite expected utility (Definition~\ref{def: completion}). 

Immediate consequences of this new result are discussed in 
Sections~\ref{sect: effective completion}--\ref{sect: intro: deflator consequences}. What emerges 
is that effective completions are linked to `corner solutions' in
the primal problem whereby, based on marginal utility considerations, the
economic agent would like to increase her exposure to risky assets in a
particular direction but cannot do so because any further exposure
takes the agent out of the effective domain of expected utility.

Our third contribution is a new construction of the optimal trading strategy
(Propositions~\ref{prop: U(X_n) converges to U(hatX) in L1(P)}, 
\ref{prop: y + Y}, \ref{prop: strategy H}, Definition~\ref{def: admissible}, and Theorem~\ref{thm: main})
where we avoid reliance on the dual optimizer altogether. This permits, for the first time
in a semimartingale setting, the construction of optimal portfolios for monotone mean-variance 
preferences. The new construction also covers the previously
unresolved case where the utility function is finite on the whole real line
but the dual optimizer is only an effective completion. We establish
existence of an optimal trading strategy under mild assumptions that reduce  
to the minimal assumptions of \citet{kramkov.schachermayer.03} 
in the $L^\infty$ case. The challenges of this construction are summarized 
in Section~\ref{sect: intro optimal strategy}.

The paper is organized as follows: in the remainder of Section~\ref{sect: intro}
we introduce the necessary concepts and notation, and discuss their economic 
and mathematical significance. Without going into too much technical detail, we also set out
the different elements of our research strategy and explain how they fit
together. Sections~\ref{sect: Fenchel duality}-\ref{sect: optimal strategy}
implement our research program and Section~\ref{sect: conclusions}
concludes. For reader's convenience Appendix~\ref{sect: duality summary}
collects useful known results in convex analysis. Appendix~\ref{sect: appxB}
constructs an explicit example of a corner solution in a continuous model
with L\'{e}vy dynamics and proves the dual optimizer cannot be linked to a
supermartingale deflator in this case. Appendix~\ref{sect: utility increases
C to C**} provides an explicit example where the duality over full market
completions fails and links it to the structure of the underlying Orlicz
space.

\subsection{Utility function \texorpdfstring{$U$}{\textit{U}} and the Orlicz
space \texorpdfstring{$L^{\hat{U}}$}{L\^{} hatU}}

\label{sect: intro utility}A utility function $U$ in this paper is a proper,
concave, non-decreasing, upper semi-continuous function. Its effective
domain is the non-empty set 
\begin{equation}\label{eq: effective domain}
\mathrm{dom}\,U=\{x\mid U(x)>-\infty \}.  
\end{equation}
The lower bound of the effective domain of $U$ is denoted by 
\begin{equation}\label{eq: underline x}
\underline{x}=\inf (\mathrm{dom}\,U).
\end{equation} 
Upper semicontinuity of $U$ means that at $\underline{x}$, which is the only possible point of discontinuity for $U$,
the utility function must be right-continuous.

The bliss point of utility is defined by 
\begin{equation}\label{eq: bliss point}
\overline{x}=\inf \{x\mid U(x)=U(\infty )\},  
\end{equation}
where $U(\infty )=\lim_{x\rightarrow \infty }U(x)$. For strictly
increasing utility functions $\overline{x}=\infty $, while for truncated
utility functions, which feature for example in monotone mean-variance
portfolio allocation \citep{cerny.al.12}, $\overline{x}<\infty $ represents
a point where further increase in wealth does not produce additional
enjoyment in terms of utility. In economics this is interpreted as the point
of maximum satisfaction, or bliss.

By construction $\underline{x}\leq \overline{x}$ and the equality arises
only when $U$ is constant on its entire effective domain in which case the
utility maximization problem is trivial because `doing nothing'\ is always
optimal. Therefore, up to a translation, the following convention entails no
loss of generality and simply means that initial endowment has been
normalized to $0$.

\begin{convention}
\label{conv: U(0)=0}{$\underline{x}<0<\overline{x}$} and $U(0)=0$.
\end{convention}

Fixing a filtered probability space $(\Omega ,\mathcal{F}_{T},P)$, the 
\emph{left} tail of the utility function $U$ gives rise to the Orlicz space 
of random variables 
\begin{equation*}
L^{\hat{U}}(\Omega ,\mathcal{F}_{T},P)=\{X\in L^{0}(\Omega ,\mathcal{F}
_{T},P)\mid E[\hat{U}\left( \lambda  X \right) ]<\infty \text{ for
some }\lambda > 0\}.
\end{equation*}
In the theory of Orlicz spaces\footnote{For a minimal overview of Orlicz 
spaces in the context of utility maximization see, for example, 
\citet[Section 2.2]{biagini.cerny.11}. A compact exposition (35 pages) appears in 
\citet[Sections 2.1 and
2.2]{edgar.sucheston.92}. Monographic references include \cite{krasnoselskii.rutickii.61} 
and \cite{rao.ren.91}.} the convex function $\hat{U}(x)=-U(-\vert x\vert ) $ is known 
as the Young function. We write $L^{\hat{U}}(P)$ or $L^{\hat{U}} $ for short when no confusion can arise.

With $X$ interpreted as the net trading gain one has $X\in L^{\hat{U}}$ if and
only if any \emph{sufficiently small} position in $X$, both long and short,
has finite expected utility. The Orlicz space $L^{\hat{U}}$ contains a smaller
subspace $M^{\hat{U}}$ (known as the Orlicz heart)\footnote{The terminology 
`Orlicz heart' appears to originate with \citet{edgar.sucheston.89}. It emphasizes 
$M^{\hat{U}}$ as a subspace of $L^{\hat{U}},$ which is a point of view important 
in our context. $M^{\hat{U}}$ can also be understood as a self-standing Banach 
space, going back to \citet[Section 8]{morse.transue.50}. Some authors use `Morse-Transue
(sub)space' or merely `Morse subspace' when referring to $M^{\hat{U}}$.} of
financial positions whose expected utility remains finite with arbitrary
scaling, 
\begin{equation*}
M^{\hat{U}}=\{X\in L^{\hat{U}}\mid E[\hat{U}\left( \lambda  X
\right) ]<\infty \text{ for all }\lambda > 0\}.
\end{equation*}

It is convenient to equip $L^{\hat{U}}$ with a Minkowski gauge norm,
\begin{equation*}
\Vert X\Vert _{\hat{U}}=\inf \{\lambda >0\mid E[\hat{U}\left( 
X /\lambda \right) ]\leq 1\},
\end{equation*}
which coincides with the classical $L^{p}$ norm when 
$\hat{U}(x )=\left\vert x\right\vert ^{p}$. In this construction the space 
$L^{\hat{U}}$ always satisfies the embeddings 
\begin{equation}
L^{\infty }\hookrightarrow L^{\hat{U}}\hookrightarrow L^{1},
\label{eq: embedding}
\end{equation}
and for quadratic utility, in particular, one obtains the natural setting
where $L^{\hat{U}}$ is isomorphic to $L^{2}$ ($L^{\hat{U}}\sim L^{2}$).

While the construction involving the space $L^{\hat{U}}$ allows one to formulate
a unified treatment for all utility functions, for topological reasons it is
at times necessary to distinguish among three cases based on the behaviour
of $U$ at $-\infty $. We flag up the three cases here for reader's
convenience.

\begin{description}
\item[Case L-F (linear, therefore finite)] 
\defcitealias{domar.musgrave.44}{Domar-Musgrave} Utility decays
asymptotically linearly, i.e. 
\begin{equation*}
0<\lim_{x\rightarrow -\infty }U(x)/x=\lim_{x\rightarrow -\infty
}U_{+}^{\prime }(x)<\infty .
\end{equation*}
A typical example is the \citetalias{domar.musgrave.44} piecewise linear
utility \citep[Fig. 3]{richter.60}. The relevant space is $L^1$.

\item[Case SL-F (super-linear and finite)] Examples include
exponential utility and truncated quadratic utility. The relevant space 
$L^{\hat{U}}$ depends on the specific $U$ but it is always strictly larger
than $L^{\infty }$ and strictly smaller than $L^{1}$,
$$L^{\infty }\hookrightdoublearrow L^{\hat{U}}\hookrightdoublearrow L^{1}.$$

\item[Case SL-INF (left tail of $U$ equals $-\infty$)] Utility
func\-tions in this category include logarithmic utility as well as power
utility functions with negative exponent. The relevant space is $L^\infty$.
\end{description}

Coarser classifications, such as \textbf{F} vs. \textbf{INF} or \textbf{L}
vs. \textbf{SL}, will be used in appropriate places. The intermediate case 
\textbf{SL}-\textbf{F} will lead to further sub-classification which will
emerge partly in the introduction and fully in the main body of the paper.
Speaking very roughly, the case $M^{\hat{U}}=L^{\hat{U}}$ will require less
work than the case $M^{\hat{U}}\subsetneq L^{\hat{U}}$; see also 
Table~\ref{tab: assumptions}.

\subsection{Primal problem and tame strategies}

The pioneering work of \cite{merton.69,merton.71,merton.73} emphasized
tractability of optimal portfolio allocation for diffusive
models of asset prices. However, \citet[Section 6]{harrison.kreps.79}
pointed out that unrestricted stochastic integration, implicit in Merton's
work, allows for so-called doubling strategies\ that lead to arbitrage
opportunities in essentially any continuous-time model of asset prices. To
prevent such economically anomalous but mathematically plausible behaviour a
consensus emerged to define `tame' strategies $\scr{T}$ as those whose wealth 
is bounded below by an arbitrary constant, see \citet[p.~400]{harrison.kreps.79}, 
\citet[Section~3.3]{harrison.pliska.81}, \citet[Theorem~1]{dybvig.huang.88}
and \citet[Definition~2.4]{karatzas.shreve.98}. We will subsequently refer to these 
strategies as $L^{\infty }$-tame,
\begin{equation*}
\scr{T}_{\infty }=\{H\in L(S)\mid \inf\nolimits_{t\in [ 0,T]}H\cdot S_{t}\in
L^{\infty }\},
\end{equation*}
where $L(S)$ is the set of all predictable $S$-integrable processes and
the symbol $H\cdot S_{t}$ stands for a stochastic integral $\int_{(0,t]}HdS$.

For a general utility function $U$ the primal portfolio allocation problem
is to compute the supremum, denoted by $u$, of expected utility over the set
of, as yet unspecified, tame trading strategies $\scr{T}$,
\begin{equation}
u(B)=\sup_{H\in \scr{T}}E[U(B+H\cdot S_{T})].  \label{eq: max over tame}
\end{equation}
Here $B\in L^{\hat{U}}$ is a random variable representing a random endowment
available at time $T$.

In the remainder of the paper $B$ is fixed and in this introduction we take 
$B=0$ for simplicity, resuming the general case from Section~\ref{sect: Fenchel duality} onwards. We
tacitly assume, along with all related literature in this area, that there
is a risk-free asset with constant value $1$ at all times. We treat the
problem (\ref{eq: max over tame}) as given, with semimartingale $S,$ utility 
$U$, and filtered probability space 
$(\Omega ,\mathbb{F}=\{\mathcal{F}_{t}\}_{t\in [ 0,T]},P)$ supplied exogenously. 
The class of tame
strategies $\scr{T}$ will be determined in response to these three inputs, 
independently of the choice of $B$.

\subsection{Tame and admissible strategies}

It is known from the deep results of \cite{kramkov.schachermayer.99,kramkov.schachermayer.03}, 
that for utility
functions with $\underline{x}$ finite (case \textbf{INF}, including log
utility and HARA\ utilities with negative exponent treated by 
\citeauthor*{merton.71}
in a lognormal setting) one can build a satisfactory framework for any
arbitrage-free semimartingale price process $S$ by restricting the agent to 
$L^{\infty }$-tame strategies. However, for utility functions with 
$\underline{x}=-\infty $ (case \textbf{F}), such as the quadratic or the negative
exponential utility, two difficulties arise that render the above approach
unsatisfactory.

The first problem is that a uniform bound on wealth rules out, for example,
normally distributed returns in any one-period model. 
\cite{biagini.frittelli.05,biagini.frittelli.07,biagini.frittelli.08}
remedy the situation by allowing tameness to depend on the utility function 
$U$, so that the maximal loss of \emph{all} tame strategies is controlled by 
\emph{one} exogenously chosen element of $L^{\hat{U}}$. In this paper we
allow losses to be controlled by \emph{any }element of $L^{\hat{U}}$.
This leads to a wider class of $L^{\hat{U}}$--tame strategies whose maximal
loss belongs to the Orlicz space $L^{\hat{U}}$,
\begin{equation*}
\scr{T}=\left\{H\in L(S)\mid \inf\nolimits_{t\in [ 0,T]}H\cdot S_{t}\in L^{\hat{U}}\right\}.
\end{equation*}

The second difficulty in the case \textbf{F} is that $L^{\infty }$--tame
strategies, and even $L^{\hat{U}}$--tame strategies as defined here, may not
contain the optimizer. Our task is to design a
larger class of \emph{admissible} strategies $\scr{A}$ 
(which again depends on $S$, $U$, and $P$, and this time also on $B$) that
attain the supremum $u(B)$ in (\ref{eq: max over tame}) without exceeding
it. This is done by requiring that each strategy in $\scr{A}$ is approximated 
in a natural sense by a sequence of strategies in $\scr{T}$, see Definition~\ref{def: admissible} 
below and also \citet[Definition~1.1]{biagini.cerny.11}
to whom we refer the reader for further background and references.

We have shown previously \citep[Theorem~4.10]{biagini.cerny.11}
that such definition of admissibility is satisfactory (at least in the case 
$B=0$) when the optimal solution in the dual problem, which we proceed to
describe below, is a $\sigma$--martingale measure. In this paper we take the
extra step to cover also the difficult case where the dual optimizer is not
a $\sigma$--martingale measure. We remark that in the setting of 
\citet{kramkov.schachermayer.99,kramkov.schachermayer.03}
the present framework yields $L^{\hat{U}}=L^{\infty }$, tame strategies $\scr{T}$ 
and admissible strategies $\scr{A}$ coincide and they are precisely those strategies 
whose wealth is bounded below by some constant.

\subsection{Economic duality \texorpdfstring{$(L^{\hat{U}},L^{\hat{V}})$}{(L\^{}hatU,L\^{}hatV)}}

\label{sect: intro economic duality}Duality has a venerable history in
economic literature. Classically, it describes the relationship between an
indirect utility function and an expenditure function or between a
production function and a cost function, see 
\citet{hotelling.32,shephard.53,blackorby.diewert.79,diewert.81}. 
Although not presented in this way historically, 
\citet{blume.08.convex_programming,blume.08.duality}
points out that microeconomic duality can be elegantly summarized using the
language of convex duality. We, too, use convex duality as a unifying theme
throughout the paper.

To avoid heavy notation, some symbols are overloaded as suggested in the
approach of \cite{rockafellar.74}. For concave $f$ the \emph{conjugate
function} $f^{\ast }$ is defined as the concave function $f^{\ast
}(y)=\inf_{x}\{\left\langle x,y\right\rangle -f(x)\}$, while the same symbol
for convex $f$ means the convex function $f^{\ast
}(y)=\sup_{x}\{\left\langle x,y\right\rangle -f(x)\}$. Here $\left\langle
x,y\right\rangle $ is a bilinear form defined over appropriate spaces, for
example $\left\langle x,y\right\rangle =xy$ when $x,y\in \mathbb{R}$.
Similarly, the \emph{effective domain }is defined as 
$\mathrm{dom}\,f=\{x:f(x)>-\infty \}$ for concave $f$ while for convex $f$ one has 
$\mathrm{dom}\,f=\{x:f(x)<\infty \}$; cf. equation (\ref{eq: effective
domain}). A self-contained technical exposition of convex duality and its
key results appears in Appendix \ref{sect: duality summary}. We suggest 
\citet{blume.08.convex_programming,blume.08.duality}
as an economic primer.

The duality constructs in this paper and related literature are somewhat
different from the classical microeconomic results surveyed above. The basic
idea here is to embed the incomplete financial market generated by tame
trading in $S$ into a statically complete financial market in which every
terminal wealth distribution in $L^{\hat{U}}$ is available at a known cost
at time $0$. This is very similar in spirit to one of the steps in the
construction of general equilibrium in \citet{arrow.debreu.54}. We now proceed 
with the detailed description of the dual pricing rules.

The concave utility $U$ has a concave conjugate $U^{\ast }$ and in line with
notation in \cite{kramkov.schachermayer.99,kramkov.schachermayer.03} we let 
$V=-U^{\ast }$. To represent pricing rules as random variables one must be
able to express prices, and with them also the bilinear form appearing in
convex duality, by an expectation operator,
\begin{equation}
\left\langle X,Y\right\rangle =E[XY],
\label{eq: bilinear form as expectation}
\end{equation}
which implies \cite[Theorems~132.2 and 132.4]{zaanen.83} that dual variables
will be taken from the Orlicz space $L^{\hat{V}}$ determined by the 
\emph{right} tail of function $V$. Here $\hat{V}\equiv \hat{U}^{\ast }$ is known
as the conjugate Young function.

Fix $Y\in L^{\hat{V}}$ and on $L^{\hat{U}}$ define a pricing rule $p_{Y}(X)=E[XY]$. 
Assuming $p_{Y}(1)=1$ one can interpret this pricing rule
as a risk-neutral expectation, $p_{Y}(X)=E^{Q_{Y}}[X]$ with $dQ_{Y}/dP=Y$.
Here $Q_{Y}$ describes time-$0$ prices of Arrow-Debreu securities in the
statically complete market with payoffs in $L^{\hat{U}}(\Omega ,\mathcal{F}
_{T},P)$ in the sense that for any contingency $A\in \mathcal{F}_{T}$ and
the corresponding elementary security with payoff $1_{A}$ the price of state 
$A$ is given by $Q_{Y}(A)\equiv E^{Q_{Y}}[1_{A}]$.

Such a market is arbitrage-free if and only if $Y>0$ $P$-a.s. which is the
same as saying that the measure $Q^{Y}$ is equivalent to the measure $P$. We denote
the set of all possible absolutely continuous state price measures by 
\begin{equation}\label{eq: P_Vhat}
P_{\hat{V}}=\left\{Q\mid dQ/dP\in L_{+}^{\hat{V}}\right\}.  
\end{equation}
We will observe later that each probability measure $Q\in P_{\hat{V}}$
describes the prices of Arrow-Debreu securities in a statically complete market
that is bliss-free for all utility functions with the same left tail and
arbitrary bounded right tail. The set of all equivalent (arbitrage-free)
pricing measures is denoted by 
\begin{equation*}
P_{\hat{V}}^{e}=\left\{Q\in P_{\hat{V}}\mid dQ/dP>0\ P\text{-a.s.}\right\}.
\end{equation*}

\subsection{Topology and duality: a caution}

In economics, duality is taken to mean a juxtaposition of two related
objects such as indirect utility function and expenditure function. In
mathematics, duality frequently refers to the choice of pairing between dual
spaces. We will now address aspects of duality in the latter sense. This
will lead to the introduction of three dual spaces with three corresponding
conjugation symbols,
\begin{equation*}
\begin{array}{ccc}
\star , & \circledast , & \ast
\end{array}
.
\end{equation*}

One typically thinks of the Orlicz space $L^{\hat{U}}$ as a Banach space
endowed with an appropriate norm. Let us denote the norm 
dual\footnote{The set of all linear functionals on $L^{\hat{U}}$ that 
are continuous for the norm topology on $L^{\hat{U}}$.} of $L^{\hat{U}}$ by $(L^{\hat{U}})^{\star }$. 
The chief difficulty facing us is that the finest topology on 
$L^{\hat{U}}$ compatible\footnote{Topology on $L^{\hat{U}}$ such that the 
set of linear functionals on $L^{\hat{U}}$ that are continuous in this 
topology can be identified with random
variables in $L^{\hat{V}}$.} with the economic duality 
$(L^{\hat{U}},L^{\hat{V}})$ of Section~\ref{sect: intro economic duality} may be strictly
coarser than the norm topology on $L^{\hat{U}}$. For example, in the case 
\textbf{INF} studied by \citet{kramkov.schachermayer.99,kramkov.schachermayer.03}
one has $L^{\hat{U}}=L^{\infty },L^{\hat{V}}=L^{1},$ while $(L^{\infty})^{\star }$ 
is the space of finitely additive measures that strictly
contains all linear functionals generated by random variables in $L^{1}$.

Topologies compatible with the economic duality 
$(L^{\hat{U}},L^{\hat{V}})$ can be characterised in more detail (see Appendices \ref{sect: duality
summary} and \ref{sect: topology equivalences}), but for now it suffices to
bear in mind that the norm topology of $L^{\hat{U}}$ may not be one of them. It
turns out that the norm topology is (trivially) compatible with economic
duality on $L^{\hat{U}}$ in finite-dimensional models while in all other
cases this happens if and only if $L^{\hat{U}}$ coincides with its Orlicz
heart $M^{\hat{U}}$ (Theorem~\ref{thm: topology equivalences}),
requiring that all financial positions in $L^{\hat{U}}$ can be scaled up
arbitrarily, long and short, while retaining finite expected utility, which
in particular covers the case \textbf{L}. It is also known that in the case 
\textbf{F }one has $(M^{\hat{U}})^{\star }=L^{\hat{V}}$ 
\citep[Theorem~2.2.11]{edgar.sucheston.92}, therefore the norm topology of 
$L^{\hat{U}}$ is compatible with economic duality on the smaller space 
$M^{\hat{U}}$, whether or not the equality $M^{\hat{U}}=L^{\hat{U}}$ holds.

The flipside is that on $L^{\hat{U}}\supsetneq M^{\hat{U}}$ one generally
loses access to many helpful properties associated with the Orlicz space 
$L^{\hat{U}}$ as a Banach space when working in the economic duality 
$(L^{\hat{U}},L^{\hat{V}})$. It also means (bi)conjugates computed in the norm topology
are in general different from those computed in the economic duality\ $(L^{\hat{U}},L^{\hat{V}})$. 
To avoid possible ambiguity we reserve the symbol $\star $ for the former and $\circledast $ for the latter. 
One should bear in mind that the bilinear form in the $\star $ duality generally cannot be
expressed by means of an expectation under measure $P$ as in equation 
(\ref{eq: bilinear form as expectation}), but merely as an abstract action of a
linear functional from $(L^{\hat{U}})^{\star }$ on an element in $L^{\hat{U}}$.

There is another level of subtlety in the general case that is not visible
in the \citet{kramkov.schachermayer.99,kramkov.schachermayer.03}
setting where $L^{\hat{U}}$ is the norm-dual of $L^{\hat{V}}$. In general
this is not true, but in the case \textbf{SL} one can recover this
relationship when $L^{\hat{V}}$ is replaced with the Orlicz heart $M^{\hat{V}}$. 
Conjugation in the duality $(L^{\hat{U}},M^{\hat{V}})$ will be denoted
by asterisk $\ast $. For readers familiar with functional analysis the
duality $(L^{\hat{U}},M^{\hat{V}})$ is compatible with the weak-star
topology on $L^{\hat{U}}.$ To summarize, in the order in which the relevant
dual spaces range from the largest to the smallest $(L^{\hat{U}})^{\star
}\hookleftarrow L^{\hat{V}}\hookleftarrow M^{\hat{V}}$ the conjugation
symbols read $\star $, $\circledast $ and $\ast $. The last two dualities
use the same bilinear form (\ref{eq: bilinear form as expectation}), given
by an expectation operator, and for any function $f$ on $L^{\hat{U}}$ the
conjugates $f^{\circledast }$ and $f^{\ast }$ coincide on $M^{\hat{V}}$.

It is known that the closure of a convex set does not depend on the choice
of a specific compatible topology, but only on the dual pair 
\citep[Theorem~5.98]{aliprantis.border.06}. To emphasize this fact we use 
the notation $\mathrm{cl}^{\star }\mathcal{A}$, $\mathrm{cl}^{\circledast }\mathcal{A}$, $\mathrm{cl}^{\ast }\mathcal{A}$, 
to denote the closure of a convex
set $\mathcal{A}$ in the three dualities. In particular, for a convex cone $\mathcal{C}
$ one has $\mathrm{cl}^{\star }\mathcal{C}=\mathcal{C}^{\star \star }$, 
$\mathrm{cl}^{\circledast }\mathcal{C}=\mathcal{C}^{\circledast\circledast }$ 
and $\mathrm{cl}^{\ast }\mathcal{C}=\mathcal{C}^{\ast \ast }$, where 
$\mathcal{C}^{\star },\mathcal{C}^{\circledast },\mathcal{C}^{\ast}$ are 
the polar cones to $\mathcal{C}$ in the appropriate duality (see equation \ref{eq: polar}).

The existing literature has very little to say about the case 
$L^{\hat{V}}\supsetneq M^{\hat{V}}$, even though this case is logically no less
important and no less prevalent in the universe of possible utility
functions than $L^{\hat{V}}=M^{\hat{V}}$. The reader is likely to be
familiar with $L^{q}$ spaces, $1\leq q<\infty $, for which one always has 
$L^{q}=M^{q}$. One might therefore think that $L^{\hat{V}}\supsetneq M^{\hat{V}}$ 
only occurs when $\hat{V}(y)$ grows faster than any power $y^{q}$.
Indeed, the right-to-left implication always holds 
\citep[eq. (I.4.7)]{krasnoselskii.rutickii.61}. But the left-to-right implication 
is not true at all; for any $q\geq 1$ one may construct $L^{\hat{V}}$ that is not 
equal to its Orlicz heart and such that 
$L^{q+\varepsilon }\hookrightarrow L^{\hat{V}}\hookrightarrow L^{q} $ 
\citep[Teorema 4]{salekhov.68}. This means that the cases $L^{\hat{V}}\supsetneq M^{\hat{V}}$ 
are interspersed in between $L^{q}$ spaces where $L^{\hat{V}}=M^{\hat{V}}$ and a theory able
to cover both is essential.

The contribution of our paper is significant already in the case $\circledast =\ast $. 
Equating $\circledast $ with $\ast $ throughout the
paper amounts to an additional assumption $L^{\hat{V}}=M^{\hat{V}}$, which
is certainly justified for all utility functions in the HARA class.

\subsection{Market completion - first attempt}

\label{sect: intro completion 1}We are now in a position to describe what,
in the field of financial economics, is classically meant by a market
completion. Denote by $\mathcal{K}$ the cone of tame terminal wealths with
zero initial capital,
\begin{equation}
\mathcal{K}=\left\{ H\cdot S_{T}:H\in \scr{T}\right\} ,  \label{eq: K}
\end{equation}
and let $\mathcal{C}$ be the convex cone of terminal wealths that are
super-replicable with zero initial capital,
\begin{equation*}
\mathcal{C}=\left( \mathcal{K}-L_{+}^{0}\right) \cap L^{\hat{U}}.
\end{equation*}

Recall the notion of Arrow-Debreu state price measure $Q\in P_{\hat{V}}$
introduced in Section~\ref{sect: intro economic duality}. With each 
$Q\in P_{\hat{V}}$, too, we associate a cone of claims that are super-replicable with
zero initial capital in the statically complete market $Q,$ 
\begin{equation}
\mathcal{C}_{Q}=\{X\in L^{\hat{U}}\mid E^{Q}[X]\equiv \left\langle
X,dQ/dP\right\rangle \leq 0\}.  \label{eq: C_Q}
\end{equation}

\begin{definition}
\label{def: separating measure}We say that the probability measure $Q$ is a
(static) \emph{market completion / separating 
measure\footnote{In the context of arbitrage theory `separation' refers to the separation of
the set of attainable claims $\mathcal{K}$ from the set of arbitrage opportunities $L_{+}^{\hat{U}}$. 
In the context of utility theory one is separating $\mathcal{K}$
from sufficiently high upper level sets of expected utility.}} if 
$\mathcal{C}\subseteq \mathcal{C}_{Q}$, which is equivalent to $dQ/dP\in \mathcal{C}^{\circledast }$, 
where $\mathcal{C}^{\circledast }$ is the polar set to $\mathcal{C}$ in the economic duality $(L^{\hat{U}},L^{\hat{V}})$,
\begin{equation}
\mathcal{C}^{\circledast }=\left\{Y\in L^{\hat{V}}\mid E[XY]\equiv \left\langle
X,Y\right\rangle \leq 0\text{ for all }X\in \mathcal{C}\subset L^{\hat{U}}\right\}.
\label{eq: polar}
\end{equation}
\end{definition}

At this point we digress a little to clarify the terminology. There is a
subtle distinction between a market \emph{completion} and a market \emph{extension} 
which disappears when $\mathcal{K}$ is a linear subspace of $L^{\hat{U}}$. 
This observation applies to any set $\mathcal{K}$ of `attainable
claims', not just the specific set in (\ref{eq: K}). It is commonly said 
\citep{ross.78,harrison.kreps.79}
that the measure $Q$ is a market extension if it correctly prices all attainable
claims, that is if $E^{Q}[X]=0$ for all $X\in \mathcal{K}$. One easily
verifies that when $\mathcal{K}$ is linear, as in the two references above,
every market completion is also a market extension and vice versa.

Earlier literature worked exclusively with linear $\mathcal{K}$. In the
context of quadratic preferences, $L^{\hat{U}}\sim L^{2}$, 
\citet{chamberlain.rothschild.83} and \citet{magill.quinzii.00}
identify the importance of continuous extension of the pricing functional to 
$\mathcal{K}^{\circledast \circledast }$. In no-arbitrage pricing, $L^{\hat{U}}\sim L^{p},p\in [ 1,\infty ]$, 
this theme is followed up by extensions to the whole of $L^{\hat{U}}$ in \citet{kreps.81}, \citet{clark.93}
and \citet{schachermayer.92,schachermayer.94}.

Starting with \citet{delbaen.schachermayer.94} and \citet{kabanov.97}
the no-arbitrage literature considers $\mathcal{K}$ that is a cone, but no
longer necessarily a linear subspace, or indeed a subset of $L^{\hat{U}}$.
This leads to situations where one may have $X\in \mathcal{K}$ such that $-X\notin \mathcal{K}$, 
as in the case of shortselling constraints\footnote{With continuous trading $\mathcal{K}$ in (\ref{eq: K}) 
may not be a linear subspace even though no explicit short-selling constraints have been imposed. This is
the situation encountered in \citet{delbaen.schachermayer.94}. 
\citet{kabanov.97} observes that one may add explicit constraints and relax 
assumption on $S$
without affecting the conclusion that $\mathcal{C}^{\star \star }\cap L_{+}^{\hat{U}}=\{0\}$ implies 
$\mathcal{C}=\mathcal{C}^{\star \star }=\mathcal{C}^{\circledast \circledast }$ when $L^{\hat{U}}\sim L^{\infty }$.}. 
To such an $X$ a separating measure may assign a strictly
negative price, $E^{Q}[X]<0$ and therefore one cannot say that $Q$ is a
`pricing measure'\ or a market extension. However, we may say that $Q$ is a
market completion because claims in $\mathcal{K}$ are attainable in the
completed market at a cost not exceeding $0$.

The mathematical necessity of using the set of super-replicable wealths $\mathcal{C}$ instead of tame 
wealths $\mathcal{K}$ stems from the fact that
one may perversely have no arbitrage over $\mathcal{K}^{\circledast
\circledast }$ while there is arbitrage over $\mathcal{C}^{\circledast
\circledast }$ \cite[Example 3.1]{schachermayer.94}. It is the statement of
the Kreps-Yan theorem\footnote{See \citet[Proposition~3.5]{gao.xanthos.17} for the Orlicz space version of
the theorem and \citet{schachermayer.02} for historical notes.} that an
arbitrage-free market completion exists ($\mathcal{C}^{\circledast }$
contains a strictly positive element) if and only if there is no arbitrage
opportunity in $\mathcal{C}^{\circledast\circledast}$ 
($\mathcal{C}^{\circledast \circledast }\cap L_{+}^{\hat{U}}=\{0\}$).

It is the statement of the even deeper Fundamental Theorem of Asset Pricing 
\citep[Theorems~1.1 and 4.1]{delbaen.schachermayer.98} that in the case $L^{\hat{U}}=L^{\infty }$ 
there is no arbitrage in $\mathcal{C}^{\circledast
\circledast }$ if and only if there is no arbitrage in the smaller
norm-closure $\mathcal{C}^{\star \star }$ and in such case 
$\mathcal{C}^{\circledast \circledast }=\mathcal{C}^{\star \star}=\mathcal{C}$ (!) and
an equivalent $\sigma $-martingale measure for $S$ exists.

\subsection{Complete market duality}

There is a dual formula \citep[Lemma~4.3]{biagini.cerny.11} that describes
the maximal utility in a complete market $Q$ in terms of its state price
density, 
\begin{equation}
u_{Q}(x)\equiv \sup_{X\in \mathcal{C}_{Q}}I_{U}(x+X)=
\min_{ \lambda\geq 0}\left\{I_{V}\left( \lambda \frac{dQ}{dP}\right)+\lambda x\right\},\text{ when }u_{Q}(x)\in\mathbb{R}.
\label{eq: u_Q(x) Fenchel}
\end{equation}
Here $I_{f}$ denotes an integral functional $I_{f}(X)=E[f(X)]$. Formula \eqref{eq: u_Q(x) Fenchel}
 arises naturally if one considers the maximization of $I_{U}(X)$ subject to a budget constraint $E^{Q}[X]=0$ 
with $\lambda $ being the Lagrange multiplier, see \citet{pliska.86}.

We say that the statically complete market $Q\in P_{\hat{V}}$ is \emph{bliss-free }if $u_{Q}(0)<U\left( \infty \right)$. 
In this case the dual formula reads 
\begin{equation}
u_{Q}(0)=\min_{\lambda >0}I_{V}\left( \lambda dQ/dP\right) .
\label{eq: u_Q Fenchel}
\end{equation}
We denote the set of all bliss-free state price measures for utility $U$ by 
\begin{equation}
P_{V}=\left\{Q\ll P\mid u_{Q}(0)<U\left( \infty \right) \right\}
=\left\{Q\ll P\mid\exists\lambda >0;I_{V}\left(\lambda \frac{dQ}{dP}\right)<\infty \right\},  \label{eq: P_V}
\end{equation}
where the set equality is hinted at in the dual formula (\ref{eq: u_Q
Fenchel}) and follows rigorously from \citet[Proposition~4.6]{biagini.cerny.11}. 
In parallel, recall the set of all absolutely continuous state price
measures $P_{\hat{V}}$ in equation (\ref{eq: P_Vhat}) and note that the
definition of the Orlicz space $L^{\hat{V}}$ allows it to be restated as 
\begin{equation}
P_{\hat{V}}=\{Q\ll P\mid \exists \lambda >0;I_{\hat{V}}(\lambda
dQ/dP)<\infty \}.  \label{eq: P_Vhat2}
\end{equation}

On comparing (\ref{eq: P_V}) and (\ref{eq: P_Vhat2}) one observes that not
all complete markets $Q\in P_{\hat{V}}$ are bliss-free because $V$ may be
unbounded near zero. It can be shown, however, that any $Q\in P_{\hat{V}}$
is bliss-free as long as $U(\infty )\equiv V(0)$ is finite, \emph{ibid}
proof i) $\Rightarrow $ ii). This underscores the economic significance of
the space $L^{\hat{V}}$ as the space of complete market pricing functionals
that are bliss-free for all utility functions sharing the same left tail and
having an arbitrary but bounded right tail. This is true for \emph{any} 
initial wealth level, as long as the initial wealth level is in the interior
of $\mathrm{dom}\,U$ and below its bliss point $\overline{x}$, \emph{ibid}.

\subsection{Martingale measures and supermartingale deflators}

\label{sect: intro deflators}So far we have suppressed the dynamic nature of
portfolio selection. To capture the temporal dimension of the problem the
no-arbitrage literature operates with $\sigma $\emph{-martingale measures}
\footnote{See \citet{emery.80} and \citet[Propositions~2.5 and 2.6]{delbaen.schachermayer.98}} 
for $S$, whose totality is denoted by 
\begin{equation}
\mathcal{M}=\{Q\ll P\mid S\text{ is a }Q\text{-}\sigma \text{--martingale}\}.  \label{eq: sigma-mart measures}
\end{equation}
Note that $S$ itself may not be a tame wealth process, that is $H=1$ may not
be a tame strategy in general. For this reason we also introduce the set of 
\emph{supermartingale measures}\footnote{Despite their superficial similarity the two notions 
`$\sigma $-martingale measures' and `supermartingale measures' refer to two very different sets of
test processes. The former relates to $S$ only; the latter refers to all
tame wealth processes $\{H\cdot S\mid H\in \scr{T}\}$.} for tame wealth processes,
\begin{equation}
\mathcal{S}=\{Q\ll P\mid H\cdot S\text{ is a }Q\text{--supermartingale for
all }H\in \scr{T}\}.  \label{eq: supermart measures}
\end{equation}

On a filtered probability space every probability measure generates
so-called \emph{density process} $\xi ^{Q}$ whose values satisfy $\xi
_{t}^{Q}=E[dQ/dP$ $|$ $\mathcal{F}_{t}]$ and therefore $\xi ^{Q}$ is a
uniformly integrable $P$-martingale. In probabilistic terms $\xi _{t}^{Q}$
is the Radon-Nikodym derivative of $Q$ restricted to $\mathcal{F}_{t}$ with
respect to $P$ restricted to $\mathcal{F}_{t}$. For an equivalent measure $Q\sim P$ 
one can use $\xi ^{Q}$ to evaluate a conditional 
price\footnote{Existence of conditional pricing rules is discussed, for example, in 
\citet{hansen.richard.87}.} $p_{t}^{Q}$ of an Arrow-Debreu security $1_{A}$ via the Bayes formula, 
\begin{equation}\label{eq: bayes}
p_{t}^{Q}(1_{A})=Q(A|\mathcal{F}_{t})\equiv E^{Q}[1_{A}|\mathcal{F}_{t}]=E[\xi _{T}^{Q}1_{A}|\mathcal{F}_{t}]/\xi _{t}^{Q}.  
\end{equation}
It follows from (\ref{eq: bayes}) that $\xi ^{Q}p^{Q}(1_{A})$ is a uniformly
integrable $P$-martingale. In these circumstances we say that $\xi ^{Q}$ is
a \emph{martingale deflator} for the price process $p^{Q}(1_{A})$. Similar
notion can be applied to the wealth of tame trading strategies.

\begin{definition}
\label{def: supermartingale deflator}Semimartingale $\xi $ is a (strong
super)martingale deflator if 
$$\text{$\xi (x+H\cdot S)$ is a $P$--(super)martingale for all $H\in \scr{T}$ and all $x\in \mathbb{R}$.}$$
 
We say that $\xi $ is a weak supermartingale deflator if instead for 
all $x\in \mathbb{R}$ the supermartingale condition holds only for some $x>0$.
\end{definition}

\begin{remark}
The set $\scr{Y}$ in \citet{kramkov.schachermayer.99,kramkov.schachermayer.03}
corresponds to the set of all supermartingale deflators for $L^{\infty }$-tame 
strategies. We will see in Section~\ref{sect: intro: deflator consequences}
that weak supermartingale deflators
are not a robust concept and only strong supermartingale deflators survive
the generalization from $L^{\infty }$ to $L^{\hat{U}}$.
\end{remark}

It turns out that each $Q\in P_{\hat{V}}\cap \mathcal{M}$ is a
supermartingale measure, $P_{\hat{V}}\cap \mathcal{M}\subseteq \mathcal{S}$
(Proposition~\ref{prop: tame strategies are super}). This in turn implies
that every $Q\in P_{\hat{V}}\cap \mathcal{M}$ is a market completion /
separating measure as per Definition~\ref{def: separating measure}, and the
cone generated by $\sigma $-martingale densities in $P_{\hat{V}}$ (denoted
with a slight abuse by $\mathcal{C}_{\sigma }^{\circledast }$),
\begin{equation}
\mathcal{C}_{\sigma }^{\circledast }=\{\lambda dQ/dP\mid \lambda \geq
0,Q\in \mathcal{M}\cap P_{\hat{V}}\},  \label{eq: sigma sep cone}
\end{equation}
is a subset of the cone $\mathcal{C}^{\circledast }$ generated by
separating densities. In the case \textbf{INF }it is additionally known that
every equivalent separating measure is a supermartingale measure 
\citep[Proposition~4.7]{delbaen.schachermayer.98}.

The converse is not true -- not every element of $\mathcal{C}^{\circledast
} $ gives rise to a $\sigma $-martingale measure for $S$ unless $S$ is
sufficiently well behaved. For our purposes it is enough to know that $\sigma $-martingale 
densities are $L^{\hat{V}}$ norm-dense in the set of separating densities. 
This is true in the case \textbf{INF} ($L^{\hat{V}}\sim L^{1}$) by \citet[Theorem~2]{kabanov.97}; 
at present the status of this conjecture in the case \textbf{F} is unknown. Therefore we make the following
\begin{assumption}
\label{ass: sig martingale measures are dense}$\mathcal{C}_{\sigma
}^{\circledast }$ is $\Vert\cdot\Vert _{\hat{V}}$-dense in $\mathcal{C}^{\circledast }$, 
that is for $Q\ll P$ with density $dQ/dP\in\mathcal{C}^{\circledast }$ 
and for every $\varepsilon >0$ there is a $\sigma $-martingale measure 
$\tilde{Q}\sim Q$ such that $\Vert d\tilde{Q}/dP-dQ/dP\Vert _{\hat{V}}\leq \varepsilon $.
\end{assumption}

There is a mild sufficient condition to guarantee that \emph{every}
separating measure is a $\sigma $-martingale measure which in turn implies
validity of Assumption~\ref{ass: sig martingale measures are dense}. For
this to hold the asset price process $S$ must be sufficiently integrable
with respect to the utility function, namely
\begin{equation}
S\in \scr{S}_{\sigma }^{\hat{U}},  \label{eq: S is nice}
\end{equation}
that is $S$ belongs $\sigma $-locally\footnote{See \citet{kallsen.04} for 
definition and properties of $\sigma $-localization. 
Further relevant properties can be found in \citet[Section~2.4]{biagini.cerny.11}.} 
to the class of processes whose maximal process at the terminal date is
in $L^{\hat{U}}$, see Sections~2.3, 2.4, Assumption~3.1, and Lemma~6.4 in 
\cite{biagini.cerny.11}\footnote{The requirement $Q\in P_{V}$ therein can be relaxed to $Q\in P_{\hat{V}}$.}.
In particular, any continuous $S$ is locally bounded ($S\in \scr{S}_{\mathrm{loc}}^{\infty }$) 
and therefore satisfies Assumption~\ref{ass: sig martingale measures are dense} for \emph{any }utility $U$ due to the
embedding (\ref{eq: embedding}).

\subsection{Duality over state price densities}\label{sect: intro Fenchel}
It is clear from the construction of a market
completion that $u_{Q}(0)$ will overestimate utility of tame trading in the
original market, $u_{Q}(0)\geq u(0)$, for any state price density $dQ/dP\in 
\mathcal{C}^{\circledast }$. We say that there is `no duality gap' if one
can complete the market in such a way that the increase in utility is
arbitrarily small,
\begin{equation}
u(0)=\sup_{X\in \mathcal{C}}I_{U}(X)=\inf_{Y\in \mathcal{C}^{\circledast
}}I_{V}(Y).  \label{eq: duality gap}
\end{equation}

It will transpire later that (\ref{eq: duality gap}) is crucial for the task
we have set out to accomplish -- which is to prove that admissible
strategies contain an optimizer. To emphasize convex duality we can write
the desirable property (\ref{eq: duality gap}) as
\begin{equation}
\sup_{X\in \mathcal{C}}I_{U}(X)=\inf_{Y\in \mathcal{C}^{\circledast
}}-I_{U^{\ast }}(Y),  \label{eq: duality over sets}
\end{equation}
or even more symmetrically as
\begin{equation}
\sup_{X\in L^{\hat{U}}}\left\{ I_{U}(X)-\delta _{\mathcal{C}}(X)\right\}
=\inf_{Y\in L^{\hat{V}}}-\left\{ I_{U}^{\circledast }(Y)
-\delta _{\mathcal{C}}^{\circledast }(Y)\right\} ,  \label{eq: Fenchel duality}
\end{equation}
where $\delta $ is convex set indicator function (zero on the set, $\infty $
outside) and $I_{U}^{\circledast }$, $\delta _{\mathcal{C}}^{\circledast }$
denote conjugate functions in the duality $(L^{\hat{U}},L^{\hat{V}})$. It is
the consequence of the careful choice of dual spaces that $I_{U}^{\circledast }=I_{U^{\ast }}=-I_{V}$. 
Because $\mathcal{C}$ is a cone, one easily obtains 
$\delta _{\mathcal{C}}^{\circledast }=\delta _{\mathcal{C}^{\circledast }}$ and 
this relationship shows equivalence between (\ref{eq: duality over sets}) and (\ref{eq: Fenchel duality}).

Results of the type (\ref{eq: Fenchel duality}) are known as the \emph{Fenchel duality}. 
For example, 
$I_{U}$ is norm-continuous at $0\in \mathcal{C}$ 
\citep[Proposition~16]{biagini.frittelli.08}
which allows application of the Fenchel duality in the norm topology 
\citep[Theorem~1.12]{brezis.11},
\begin{equation}\label{eq: FM strong}
\sup_{X\in L^{\hat{U}}}\left\{ I_{U}(X)-\delta _{\mathcal{C}}(X)\right\}
=\min_{Y\in (L^{\hat{U}})^{\star }}-\left\{ I_{U}^{\star }(Y)-\delta _{\mathcal{C}}^{\star }(Y)\right\} .  
\end{equation}
This is formally the same formula as (\ref{eq: Fenchel duality}) but with a
larger dual space.

In our previous work we had to assume that the dual minimizer on the
right-hand side of (\ref{eq: FM strong}) was an element of $\mathcal{C}^{\circledast }$, 
i.e. a separating measure. Here we remove that
assumption. Our first step is to rewrite the known result (\ref{eq: FM
strong}) for the norm duality $(L^{\hat{U}},(L^{\hat{U}})^{\star })$, in
terms of the economic duality $(L^{\hat{U}},L^{\hat{V}})$, 
\begin{equation}
\sup_{X\in L^{\hat{U}}}\left\{ I_{U}(X)-\delta _{\mathcal{C}\cap \mathcal{D}}(X)\right\} 
=\min_{Y\in L_{+}^{\hat{V}}}-\left\{ I_{U}^{\circledast}(Y)-\delta _{\mathcal{C}\cap \mathcal{D}}^{\circledast }(Y)\right\} ,
\label{eq: dual corner}
\end{equation}
where $\mathcal{D}=\mathrm{dom}\,I_{U}$ is the effective domain of expected
utility (see Theorem~\ref{thm: strong duality}).

Crucially for our story $\mathcal{C}\cap \mathcal{D}$ may be a
strict subset of $\mathcal{C}$ and therefore the dual optimizer in (\ref{eq:
dual corner}) need not be an element of $\mathcal{C}^{\circledast }$ and
thus not a separating measure, even though it necessarily must be a state price density
in $L_{+}^{\hat{V}}$.

\subsection{Optimal trading strategy}\label{sect: intro optimal strategy}
Armed with the previous observation, we adopt a radically new approach that bypasses 
the dual optimizer entirely. Instead, we construct 
the candidate optimal trading strategy from 
a supermartingale compactness result of \citet[Theorem D]{delbaen.schachermayer.98}, 
using one arbitrary $\sigma $-martingale measure whose existence we assume. Having 
proved in Section~\ref{sect: optimal wealth} that the utility of wealth of the 
maximizing sequence can be chosen to have an integrable lower bound 
(Proposition~\ref{prop: U(X_n) converges to U(hatX) in L1(P)}) the difficulty is then 
showing that the expected utility of the candidate terminal wealth does not \emph{exceed} 
$u(0)$, the expected utility attainable by tame trading.

By carefully rethinking the arguments of \citet[Proposition~3.8]{biagini.cerny.11}
we observe that the candidate optimal wealth process is a supermartingale
under every $\sigma$--martingale measure. Consequently, the
utility of the candidate wealth is majorized by the utility of every $\sigma$--martingale 
measure and the new construction goes through as long as \emph{there is
no duality gap over separating measures}, that is (\ref{eq: duality gap})
holds, and $\sigma$--martingale measures are suitably dense among separating measures. 
The desired duality is proved in Section~\ref{sect: duality over separating}. 

The advantage of the proposed construction is twofold --- it allows us to 
deal with the case when the dual optimizer is not a separating measure 
and it also covers the case where the
optimal wealth is in the algebraic interior of the effective domain but the
dual maximizer, which now must be\emph{\ }a separating measure, is not
equivalent to $P$. In the latter case the utility function is not strictly
monotone. To give an example, truncated quadratic utility plays an important role
in the computation of monotone mean-variance optimal portfolios, see \cite{maccheroni.al.09} and \cite{cerny.al.12}.

The remaining sections implement the research program outlined above.

\section{Fenchel duality over state price densities}\label{sect: Fenchel duality}

To allow for random endowments, the set $\mathcal{C}$ is replaced
by the set $B+\mathcal{C}$.
Because $0$ may not be an element of $B+\mathcal{C}$, the
arguments leading to (\ref{eq: FM strong}) may fail. However, $I_{U}$ is
norm-continuous not only at zero but everywhere on the algebraic interior of 
$\mathcal{D}=\mathrm{dom}\,I_U$ \citep[Corollary 8B]{rockafellar.74}. In the present setting 
the algebraic interior is given explicitly as
\begin{equation}
\mathrm{core\,}\mathcal{D}=\{X\in L^{\hat{U}}:\exists \lambda
>1;I_{U}(\lambda X)>-\infty \}\text{.}  \label{eq: core(D)}
\end{equation}
Financially these are the positions that allow for proportional increase
while maintaining finite utility level.

\begin{theorem}
\label{thm: strong duality}
Assume $B\in L^{\hat{U}}$, $(B+\mathcal{C})\cap \mathrm{core\,}\mathcal{D\neq \emptyset }$, and let 
$\mathcal{A}=B+\mathcal{C}$. Then one has
\begin{equation}\label{eq: duality AcapD}
\begin{split}
\sup_{X\in \mathcal{A}}I_{U}(X) &=\sup_{X\in \mathrm{cl}^{\star }\mathcal{A}}I_{U}(X)
=\sup_{X\in \mathcal{A}\cap \mathcal{D}}I_{U}(X)
=\sup_{X\in \mathrm{cl}^{\circledast }\left( \mathcal{A}\cap \mathcal{D}\right) }I_{U}(X)\\
&=\min_{Y\in L_{+}^{\hat{V}}}\left\{ I_{V}\left( Y\right) +\delta_{\mathcal{A}\cap \mathcal{D}}^{\circledast }(Y)\right\},
\end{split}
\end{equation}
where the so-called support function $\delta _{\mathcal{G}}^{\circledast
}(Y)$ has the explicit form 
\begin{equation*}
\delta _{\mathcal{G}}^{\circledast }(Y)=\sup_{X\in \mathcal{G}}E\left[ XY\right] .
\end{equation*}
\end{theorem}

\begin{proof}
\emph{Step 1) }The Fenchel inequality $I_{U}(X)\leq I_{V}(Y)+E[XY]$ for $X\in L^{\hat{U}},Y\in L_{+}^{\hat{V}}$ yields 
\begin{equation}
u(B)=\sup_{X\in \mathcal{A}\cap \mathcal{D}}I_{U}(X)\leq I_{V}(Y)+\sup_{X\in 
\mathcal{A}\cap \mathcal{D}}E[XY]\text{ for all }Y\in L_{+}^{\hat{V}}.
\label{dual-2}
\end{equation}
When $u(B)=\infty $ we necessarily have $U(\infty )\equiv V(0)=\infty $ and
the duality (\ref{eq: duality AcapD}) therefore holds trivially with $Y=0$.

\emph{Step 2) }Consider the remaining case $u(B)<\infty $. The Fenchel
duality in the norm topology \citep[Theorem 1.12]{brezis.11} gives
\begin{equation}\label{eq: norm duality}
\begin{split}
u(B) &=\sup_{X\in \mathcal{A}}I_{U}(X)=\min_{\mu \in (L^{\hat{U}})^{\star
}}-\left\{ I_{U}^{\star }(\mu )-\delta _{\mathcal{A}}^{\star }(\mu )\right\}
\\
&=\min_{\mu \in (L^{\hat{U}})^{\star }}\{-I_{U}^{\star }(\mu )+\mu
(B)+\sup_{X\in \mathcal{C}}\mu (X)\},
\end{split}
\end{equation}
where $\mu (X)=\int X(\omega )\mu (d\omega )$. We now invoke finiteness of $u(B)$ and observe, 
because $\mathcal{C}$ is a cone, that the right-hand side is finite only if $\mu \in \mathcal{C}^{\star }$ which yields
\begin{equation}
u(B)=\min_{\mu \in \mathcal{C}^{\star }}\{-I_{U}^{\star }(\mu )+\mu (B)\}.
\label{dual-1}
\end{equation}
One can repeat the same argument starting with $\mathrm{cl}^{\star }\mathcal{A}=B+\mathrm{cl}^{\star }\mathcal{C}$ 
in place of $\mathcal{A}$ to find the right-hand side in (\ref{dual-1}) remains unchanged. This proves 
$$\sup_{X\in\mathcal{A}}I_{U}(X)=\sup_{X\in \mathrm{cl}^{\star }\mathcal{A}}I_{U}(X).$$

\emph{Step 3) }By \citet[Theorem~2.6]{kozek.79} the conjugate $I_{U}^{\star }$ 
on the norm-dual of $L^{\hat{U}}$ is given explicitly by
\begin{equation*}
-I_{U}^{\star }(\mu )=I_{-U^{\ast }}(d\mu _{r}/dP)+\delta _{\mathcal{D}
}^{\star }(-\mu _{s}),
\end{equation*}
where $\delta _{\mathcal{D}}^{\star }(\mu )=\sup_{X\in \mathcal{D}}\mu (X)$
is the convex conjugate of the convex indicator function $\delta _{\mathcal{D}}$ 
\citep[equation~(3.13)]{rockafellar.74}; $\mu =\mu _{r}+\mu _{s}$ is a unique decomposition 
of $\mu $ into a regular and singular part 
\citep[Theorem~133.6]{zaanen.83}; and $Y=d\mu _{r}/dP\in L^{\hat{V}}$. 
Therefore (\ref{dual-1}) can be written as
\begin{equation}
u(B)=\sup_{X\in \mathcal{A}}I_{U}(X)=\min_{\mu \in \mathcal{C}^{\star
}}\{I_{V}(Y)+\mu (B)+\delta _{\mathcal{D}}^{\star }(-\mu _{s})\},
\label{dual0}
\end{equation}
see also \citet[Theorem~3.8]{biagini.al.11}.

\emph{Step 3) }Denote the minimizer on the right-hand side of (\ref{dual0})
by $\hat{\mu},$ with $\hat{Y}=\frac{d\hat{\mu}_{r}}{dP}$,
\begin{equation}\label{dual1}
\begin{split}
u(B) &=I_{V}(\hat{Y})+\hat{\mu}(B)+\sup_{X\in \mathcal{D}}-\hat{\mu}_{s}(X)\\
&=I_{V}(\hat{Y})+\hat{\mu}(B)+\sup_{X\in \mathcal{D}-B}-\hat{\mu}_{s}(X+B) \\
&=I_{V}(\hat{Y})+E[\hat{Y}B]+\sup_{X\in \mathcal{D}-B}-\hat{\mu}_{s}(X). 
\end{split}
\end{equation}
Rephrase (\ref{dual-2}) as
\begin{equation}
u(B)=\sup_{X\in \mathcal{A}\cap \mathcal{D}-B}I_{U}(X+B)\leq I_{V}(\hat{Y})
+\sup_{X\in \mathcal{A}\cap \mathcal{D}-B}E[(X+B)\hat{Y}],  \label{dual2}
\end{equation}
and combine this with (\ref{dual1}) to obtain 
\begin{equation}
\sup_{X\in \mathcal{D}-B}-\hat{\mu}_{s}(X)\leq \sup_{X\in \mathcal{A}\cap 
\mathcal{D}-B}E[X\hat{Y}].  \label{dual3}
\end{equation}

\emph{\ Step 4) }Recall the notation $\mathcal{A}=B+\mathcal{C}.$ Recall $\hat{\mu}\in \mathcal{C}^{\star }$ 
and, because of the polar relationship $\hat{\mu}(X)\leq 0$ for all $X\in \mathcal{C}$, we have 
$-\hat{\mu}_{s}(X)\geq E[X\hat{Y}]$ for all $X\in \mathcal{C}\supset \mathcal{A}\cap\mathcal{D}-B$ which yields 
\begin{equation}
\sup_{X\in \mathcal{A}\cap \mathcal{D}-B}-\hat{\mu}_{s}(X)\geq \sup_{X\in\mathcal{A}\cap \mathcal{D}-B}E[X\hat{Y}].  \label{dual4}
\end{equation}
From (\ref{dual3}-\ref{dual4}) we obtain the following chain of inequalities
\begin{equation}\label{eq: norm of yhat_s}
\begin{split}
\sup_{X\in \mathcal{D}-B}-\hat{\mu}_{s}(X)&\geq \sup_{X\in \mathcal{A}\cap 
\mathcal{D}-B}-\hat{\mu}_{s}(X)\\
&\geq \sup_{X\in \mathcal{A}\cap \mathcal{D}-B}E[X\hat{Y}]\geq \sup_{X\in \mathcal{D}-B}-\hat{\mu}_{s}(X),
\end{split}
\end{equation}
which are therefore equalities. On combining (\ref{dual1}) and (\ref{eq:
norm of yhat_s}), together with an explicit expression for the support
function \citep[equation (3.13)]{rockafellar.74}
one obtains equality in (\ref{dual2}),
\begin{equation}
u(B)=I_{V}(\hat{Y})+E[\hat{Y}B]+\sup_{X\in \mathcal{A}\cap \mathcal{D}-B}E[X\hat{Y}]
=I_{V}(\hat{Y})+\delta _{\mathcal{A}\cap \mathcal{D}}^{\circledast}(\hat{Y}).  \label{dual5}
\end{equation}

\emph{Step 5) }By continuity of the bilinear form $\langle X,\hat{Y}\rangle
\equiv E[X\hat{Y}]$ in the $\circledast $ duality one has 
\begin{equation*}
\sup_{X\in \mathcal{A}\cap \mathcal{D}}E[X\hat{Y}]=\sup_{X\in \mathrm{cl}^{\circledast }(\mathcal{A}\cap \mathcal{D)}}E[X\hat{Y}],
\end{equation*}
which when combined with the Fenchel inequality and (\ref{dual5}) yields
\begin{equation*}
u(B)\leq \sup_{X\in \mathrm{cl}^{\circledast }(\mathcal{A}\cap \mathcal{D)}}I_{U}(X)
\leq I_{V}(\hat{Y})+\sup_{X\in \mathrm{cl}^{\circledast }(\mathcal{A}\cap \mathcal{D)}}E[X\hat{Y}]=u(B).
\end{equation*}
This completes the proof in the remaining case $u(B)<\infty$.
\end{proof}

Observe that Theorem~\ref{thm: strong duality} does not claim $\sup_{X\in 
\mathcal{A}}I_{U}(X)=\sup_{X\in \mathrm{cl}^{\circledast }\mathcal{A}}I_{U}(X)$. 
Appendix \ref{sect: utility increases} gives an example with $B=0 $ where one obtains strict inequality 
$\sup_{X\in \mathcal{C}}I_{U}(X)<\sup_{X\in \mathcal{C}^{\circledast \circledast }}I_{U}(X)$.
Nonetheless, Theorem~\ref{thm: strong duality} continues to hold for $\mathcal{A}=\mathcal{C}$ 
as well as for $\mathcal{A}=\mathcal{C}^{\circledast \circledast }$ except each case must by necessity have a
different dual optimizer.

We remark that (\ref{eq: duality AcapD}) can be equivalently rephrased as
\begin{equation}\label{eq: dual corner 2}
\sup_{X\in L^{\hat{U}}}\left\{ I_{U}(X)-\delta _{\mathcal{A}\cap \mathcal{D}}(X)\right\} 
=\min_{Y\in L^{\hat{V}}}-\left\{ I_{U}^{\circledast}(Y)-\delta _{\mathcal{A}\cap \mathcal{D}}^{\circledast }(Y)\right\},
\end{equation}
which signifies that the left-hand side and the right-hand side form a
relationship known as the strong Fenchel duality. The new result (\ref{eq: dual corner 2}) 
is mathematically significant because standard regularity
conditions for the Fenchel duality require $L^{\hat{U}}$ to be normed or at
least metric while in the strongest available topology for the pair 
$(L^{\hat{U}},L^{\hat{V}})$ the space $L^{\hat{U}}$ generally fails to be
barreled (\emph{tonnel\'{e}} in \citealp{rockafellar.66}) and therefore 
cannot be compatible with the metric or the norm topology. 

\cite{bot.10} Theorems~2.2, 
15.2, and Remark~7.8 summarize regularity conditions
under which (\ref{eq: dual corner 2}) is known to hold, but in the present
case none of these conditions applies. The conditions in Theorem~2.2 are not
applicable because $L^{\hat{U}}$ may not be a Riesz (metric) space in any
topology compatible with duality $(L^{\hat{U}},L^{\hat{V}})$; those in
Remark~7.8 and Theorem~15.2 fail because $\mathcal{C}$ is not necessarily $\circledast$-closed.

In the case \textbf{F }with $B=0$ the dual formula (\ref{eq: dual corner}) was
obtained independently by \citet{gushchin.al.14}. In comparison, our approach 
is more direct, covering both \textbf{F} and 
\textbf{INF} case in one go and producing a proof that is, even just in the 
\textbf{F} case, significantly shorter, while allowing for random endowment.

In the literature on utility maximization with random endowment the case 
\textbf{INF }is covered by \citet{cvitanic.al.01}
who assume $B\in L^{\infty }$ therefore $x+B\in \mathrm{core\,}\mathcal{D}$
for $\underline{x}<x$. Their dual formula, containing singular parts,
corresponds to our equation (\ref{dual0}) with $B$ replaced by $x+B$ once we
realize that in their setting $\mathcal{D}$ is the set of strictly positive
random variables in $L^{\infty }$ and therefore $\delta _{\mathcal{D}}^{\star }(-\mu _{s})=0$. 
In the same setting, \citet{hugonnier.kramkov.04} remove the singular parts from the dual, 
using methods similar to those of \citet{kramkov.schachermayer.03}.

In the case \textbf{F}, \citet[Definition~3.1]{biagini.al.11}
have a condition equivalent to $B\in \mathrm{core\,}\mathcal{D}$ which is
stronger than our assumption $(B+\mathcal{C})\cap \mathrm{core\,}\mathcal{D}\neq \emptyset$, 
and just like \citet{cvitanic.al.01} their dual problem contains singular parts. 

For an immediate consequence of duality (\ref{eq: dual corner 2}) recall that
the largest linear subspace of $L^{\hat{U}}$ contained in $\mathrm{dom}\,I_{\hat{U}}$ 
is known as the Orlicz heart $M^{\hat{U}}$. It is now evident from (\ref{eq: dual corner 2}) 
and from the inclusion $\mathrm{dom}\,I_{\hat{U}}\subseteq \mathrm{dom}\,I_{U}\equiv \mathcal{D}$ 
that the dual optimizer will correspond to a separating measure if $L^{\hat{U}}=M^{\hat{U}}$ or if
at least $B+\mathcal{C}\subseteq M^{\hat{U}}$ because then $(B+\mathcal{C})\cap\mathcal{D}=B+\mathcal{C}$. 
In particular, when working with locally bounded processes one may opt for $L^{\infty }$-tame strategies 
controlled from both sides \citep{biagini.cerny.11} whereby $\mathcal{C}\subseteq L^{\infty }$. 
In the case \textbf{F} one has $L^{\infty }\hookrightarrow M^{\hat{U}}$ and the
Fenchel duality (\ref{eq: dual corner 2}) then yields for any $B\in M^{\hat{U}}$
a utility-based Fundamental Theorem of Asset Pricing, previously obtained under 
an additional assumption stronger than $M^{\hat{V}}=L^{\hat{V}}$ in \citet[Theorem~1.2]{owen.zitkovic.09}.

Further important consequences of the new formula (\ref{eq: dual corner 2})
are described in the next two sections.

\subsection{Market completion - a new definition}\label{sect: effective completion}

Denoting an optimizer on the right-hand side of the Fenchel duality (\ref{eq: dual corner}) 
by $\hat{Y}$ and setting $d\hat{Q}/dP=\hat{Y}/E[\hat{Y}]$, in view of equation 
(\ref{eq: u_Q(x) Fenchel}) we may interpret the
right-hand side expression in (\ref{eq: dual corner}) as the maximal utility in
a bliss-free complete market $\hat{Q}$ with initial endowment increased by
the amount $\delta _{\mathcal{C}\cap \mathcal{D}}^{\circledast }(d\hat{Q}/dP)$. 
This market completion is somewhat unusual because we are not
completing the entire market $\mathcal{C}$, merely the part where the
expected utility is finite, $\mathcal{C}\cap \mathcal{D}$, and extra initial
endowment is required.

\begin{definition}\label{def: completion}
We say $Q\in P_{\hat{V}}$ is a completion of market $\mathcal{C}$ if 
$\mathcal{C}\cap \mathcal{D}\subseteq x+\mathcal{C}_{Q}$ for some $x\in[ 0,\infty )$. 
When $x$ can be chosen equal to zero we say $Q$ is a \emph{full completion}, 
otherwise we say $Q$ is an \emph{effective completion}.
\end{definition}

It follows that $Q$ is a completion if and only if 
$dQ/dP\in \mathrm{dom}\,\delta _{\mathcal{C}\cap \mathcal{D}}^{\circledast }\supseteq 
\mathrm{dom}\,\delta _{\mathcal{C}}^{\circledast }\equiv \mathcal{C}^{\circledast }$.
The terminology \emph{full completion} is justified by the equivalence 
\begin{equation}\label{eq: full completion equivalence}
\mathcal{C}\cap \mathcal{D}\subseteq \mathcal{C}_{Q}\Leftrightarrow \mathcal{C}\subseteq \mathcal{C}_{Q},  
\end{equation}
which follows from the observation that $0$ is in the norm interior of $\mathcal{D}$ 
implying $\mathcal{C}=\mathrm{cone}\left( \mathcal{C}\cap\mathcal{D}\right) $. 
A full completion $Q$ is therefore precisely the classical completion discussed in Section~\ref{sect: intro completion 1},
that is a separating measure.

We can now interpret the Fenchel duality formula (\ref{eq: dual corner}) as
a market completion theorem: market $\mathcal{C}$ is bliss-free if and only
if $\mathcal{C}\cap \mathcal{D}$ can be embedded in a bliss-free complete
market with the same expected utility.

\subsection{Boundary solutions, corner solutions, and separating measures}

\label{sect: intro shadow price}For the purpose of this section we assume $V$
is strictly convex on $\mathrm{dom}\,V$. Because $V$ is closed 
(Proposition~\ref{prop: h*}), it follows by \citet[Theorem~26.3]{rockafellar.70}
this is equivalent to $U$ being essentially smooth, that is differentiable
on $(\underline{x},\infty )$ and satisfying $\lim_{x\searrow \underline{x}}U^{\prime }(x)=\infty $, 
which explains the origin of technical conditions
customarily imposed on $U$ in the literature. Strict convexity of $V$ means
the dual optimizer in (\ref{eq: dual corner}) is necessarily unique. We
denote it by $\hat{Y}$ and let $$d\hat{Q}/dP=\hat{Y}/E[\hat{Y}].$$

While the emergence of the set $\mathcal{C}\cap \mathcal{D}$ in formula 
(\ref{eq: dual corner}) is unexpected, \emph{post hoc} it has a natural economic
interpretation. The fact that $\mathcal{C}\cap \mathcal{D}$ may be a strict
subset of $\mathcal{C}$ implies that the primal optimum $\hat{X}$ (supposing
it exists in $L^{1}(\hat{Q})$-closure of $\mathcal{C}\cap \mathcal{D}$ as
discussed in Section~\ref{sect: optimal wealth}) may be a `boundary
solution'\ in the sense that $\theta \hat{X}\notin \mathcal{D}$ for $\theta
>1$.

Consider the constrained optimization $\max_{\theta \leq 1}I_{U}(\theta \hat{X})$. 
When $E[\hat{X}U^{\prime }(\hat{X})]>0$ the constraint $\theta \leq 1$
is binding and the Lagrange multiplier associated with this constraint is
exactly equal to $E[\hat{X}U^{\prime }(\hat{X})]$. We will refer to this
situation as a `corner solution'. It is interesting to note that the
constraint in question is not exogenous, rather the corner arises implicitly
due to the boundedness of the effective domain $\mathcal{D}$ in some
directions. Under the current hypotheses the following statements are
equivalent:

\begin{enumerate}
\item $\hat{Q}$ is not a separating measure (i.e. $\hat{Q}$ is only an
effective completion);

\item $\delta _{\mathcal{C}\cap \mathcal{D}}^{\circledast }(\hat{Y})=E[\hat{X}\hat{Y}]=E[U^{\prime }(\hat{X})\hat{X}]>0$.
\end{enumerate}

As can be expected, boundary solution\ is \emph{not synonymous} with 
$\delta_{\mathcal{C}\cap \mathcal{D}}^{\circledast }(\hat{Y})>0$. In particular,
when optimizing over a complete market one always has 
$\delta _{\mathcal{C}\cap \mathcal{D}}^{\circledast }(\hat{Y})=0$ 
(see equation \ref{eq: u_Q Fenchel}), while the primal solution may lie on the edge of the effective
domain $\mathcal{D}$. The converse statement that non-boundary primal optimizer corresponds to a 
full completion in the dual problem appears, with an extra technical condition, in \citet[Proposition~31]{biagini.frittelli.08}.

\subsection{Implications for supermartingale deflators}\label{sect: intro: deflator consequences}
In \citet{kramkov.schachermayer.99,kramkov.schachermayer.03}
the dual optimizer $\hat{Y}$ is interpreted as a terminal value of a
supermartingale deflator. Now suppose that the optimal terminal wealth $\hat{X}$ introduced 
in Section~\ref{sect: intro shadow price} has an optimal
strategy $\hat{H}$ associated with it, $\hat{X}=\hat{H}\cdot S_{T}$. We
already know that $\hat{Y}$ is an effective completion (that is, not a
separating measure)\ if and only if $\hat{X}$ is a corner solution, in which
case
\begin{equation*}
\hat{H}\cdot S_{0}=0<E[\hat{X}\hat{Y}]=E[\hat{Y}(\hat{H}\cdot S_{T})].
\end{equation*}
This inequality means that there can be \emph{no strong }supermartingale
deflator for $\hat{H}\cdot S$ with terminal value $\hat{Y},$ if $\hat{Y}$ is
an effective completion.

In \citet{kramkov.schachermayer.99,kramkov.schachermayer.03}
(case \textbf{INF}) an effective completion $\hat{Y}$ can be turned into a 
\emph{weak }supermartingale deflator by setting 
\begin{equation*}
\xi _{t}=\frac{E[\hat{Y}(x+\hat{H}\cdot S_{T})|\mathcal{F}_{t}]}{x+\hat{H}\cdot S_{t}},
\end{equation*}
with $\underline{x}=0$ and $x>0$ in the notation of Section~\ref{sect: intro utility}. 
However, when $U$ is finite everywhere (case \textbf{F}) 
for $\hat{Y}\notin \mathcal{C}^{\circledast }$ there may be \emph{no}
supermartingale deflator with terminal value $\hat{Y}$ at all. An example of
such situation is given in Appendix \ref{sect: no weak supermartingale
deflator}. This shows that an association of a supermartingale deflator with a
dual optimizer not in $\mathcal{C}^{\circledast}$ is an ad-hoc construction.

One can robustly characterize the effective completion $\hat{Q}$ as a 
\emph{``submartingale''} measure in the sense
that for the optimal trading strategy $\hat{H}$ one will have 
\begin{equation*}
E^{\hat{Q}}[\hat{H}\cdot S_{T}]=\delta _{\mathcal{C}\cap \mathcal{D}}^{\circledast }(d\hat{Q}/dP) 
= \sup_{X\in \mathcal{C}\cap \mathcal{D}}E^{\hat{Q}}[X]>0 = \hat{H}\cdot S_{0}.
\end{equation*}
We conjecture that for the effective completion $\hat{Q}$ the submartingale
property $E^{\hat{Q}}[\hat{H}\cdot S_{u}|\mathcal{F}_{t}]\geq \hat{H}\cdot
S_{t\wedge u}$ holds for $u=T$ and arbitrary $t$ although not necessarily
for all $u$ and $t$. The relationship $E^{\hat{Q}}[\hat{H}\cdot S_{T}]>\hat{H}\cdot S_{0}$ 
is universal across utility functions and robust to arbitrary
translation of initial wealth.

\section{Optimal terminal wealth}

\label{sect: optimal wealth}

Our next step is to show that there is an optimizing sequence $\{X_{n}\}$ of
terminal wealth distributions in $\mathcal{C}$ which converges pointwise $P$-a.s. 
to a limit $\hat{X}$ and such that $U(B+X_{n})$ approximate $U(B+\hat{X})$ in $L^{1}(P)$. 
This means that $\hat{X}$ necessarily attains the maximal
utility $u(B)$. The desired convergence requires uniform integrability of
the sequence $\{U(B+X_{n})\}$ which in general fails to materialize, even in
the `nice' case $M^{\hat{U}}=L^{\hat{U}}$.

A complete study of minimal conditions for the uniform integrability of the
utility of maximizing sequence is beyond the scope of this paper. We remark
that the key tools in that direction are the results of \citet{ando.62}
on compactness in the economic duality $(L^{\hat{U}},L^{\hat{V}})$. In this
section we proceed by introducing comparatively simple sufficient conditions
encompassing all results available to date.

Denote $a=U_{+}^{\prime }(0)$. Recall that $V=-U^{\ast }$ and due to 
$U(0)=0 $ function $V$ is decreasing on $[0,a]$ and increasing on $[a,\infty
) $. Recall $\hat{U}(x)=-U\left(-\lvert x\rvert\right)$ and $X\in L^{\hat{U}}$ 
if there is $\lambda >0$ such that $I_{\hat{U}}\left(\lambda X \right) <\infty $. 
We have $\hat{V}(y)=\hat{U}^{\ast }(y)=V(\lvert y\rvert\vee a)$.

The first important ingredient is the requirement that the elements of $B+\mathcal{C}$ 
with high expected utility must have negative parts of bounded $L^{\hat{U}}$ norm. 
This requirement is satisfied trivially in the case 
\textbf{INF (}$L^{\hat{U}}\sim L^{\infty }$) and not just over $B+\mathcal{C}$ but over the 
entire space $L^{\infty }$ because in that case 
$$I_{U}(X)>-\infty \Rightarrow X\geq \underline{x}\iff \Vert X^{-}\Vert_{\infty }\leq -\underline{x}.$$ 
The following concept appears to be new.

\begin{definition}
\label{def: norm-coercive in losses}We say that expected utility $I_{U}$ is 
\emph{norm-coercive in losses on a set} $\mathcal{G}\subseteq L^{\hat{U}}$ if
\begin{equation}
\lim_{\Vert X^{-}\Vert _{\hat{U}}\rightarrow \infty ,X\in \mathcal{G}}I_{U}(X)=-\infty .  \label{eq: norm-coercive in losses}
\end{equation}
Equivalently, expected utility is norm-coercive in losses on $\mathcal{G}$
if and only if for every $k\in \mathbb{R}$ there is $l>0$ such that $I_{U}(X)>k$ 
implies $\Vert X^{-}\Vert _{\hat{U}}\leq l$ for all $X\in \mathcal{G}$.
\end{definition}

We continue with a lemma that establishes boundedness properties for
cost-constrained subsets of upper level sets of expected utility and leads
to sufficient conditions that imply norm coercivity in losses. For $U$
bounded above, expected utility is trivially norm-coercive in losses over
the entire space $L^{\hat{U}}$. This can be seen also in the lemma below by
setting $\tilde{Y}=0$, which is possible in the bounded case thanks to $U\left( \infty \right) = V(0)<\infty $.

\begin{lemma}
\label{lem: L^1(P) boundedness}Consider a set $\mathcal{G}\subseteq L^{\hat{U}}$ 
and suppose there is $\tilde{Y}\in L^{\hat{V}}$ such that $\{E[X\tilde{Y}]\}_{X\in \mathcal{G}}$ 
is bounded from above, and $\lambda \tilde{Y}\in\mathrm{dom}\,I_{V}$ for two distinct 
values of $\lambda >0$. Consider further an arbitrary $\widetilde{\mathcal{G}}\subseteq \mathcal{G}$ such
that $\{I_{U}(X)\}_{X\in \widetilde{\mathcal{G}}}$ is bounded from below.
The following statements hold:

i) $\{\left\vert X\right\vert \tilde{Y}\}_{X\in \widetilde{\mathcal{G}}}$ is 
$L^{1}(P)$-bounded;

ii) $\{\left\vert U(X)\right\vert \}_{X\in \widetilde{\mathcal{G}}}$ is $L^{1}(P)$-bounded;

iii) $I_{U}$ is norm-coercive in losses on $\mathcal{G}$.
\end{lemma}

\begin{proof}
\emph{i)} Consider $\lambda _{2}>\lambda _{1}>0$ such that $\lambda _{i}\tilde{Y}\in \mathrm{dom}\,I_{V}$, $i=1,2$. 
By the Fenchel inequality
\begin{align}
U\left( X^{+}\right) &\leq V(\lambda _{1}Y)+\lambda _{1}X^{+}Y,
\label{eq: Fenchel Uplus} \\
U\left( -X^{-}\right) &\leq V\left( \left( \lambda _{2}Y\right) \vee
a\right) -\lambda _{2}X^{-}Y.  \label{eq: Fenchel Uminus}
\end{align}
On taking expectations
\begin{equation}
I_{U}(X)\leq I_{V}(\lambda _{1}\tilde{Y})+I_{V}(\lambda _{2}\tilde{Y}
)+\lambda _{1}E[X\tilde{Y}]-\left( \lambda _{2}-\lambda _{1}\right) E[X^{-}
\tilde{Y}].  \label{eq: key inequality}
\end{equation}
Note that $\tilde{Y}\in L^{\hat{V}}$ and \ $|X\tilde{Y}|\in L^{1}(P)$ for
any $X\in L^{\hat{U}}$ by the Orlicz space H\"{o}lder inequality 
\citep[eq.~3.3.4]{rao.ren.91}. As $\{I_{U}(X)\}_{X\in \widetilde{\mathcal{G}}}$ is bounded below, the
assumed upper bound on $\{E[X\tilde{Y}]\}_{X\in \mathcal{G}}$ and (\ref{eq:
key inequality}) imply $\{X^{-}\tilde{Y}\}_{X\in \widetilde{\mathcal{G}}}$
is $L^{1}(P)$-bounded, therefore $\{X^{+}\tilde{Y}\}_{X\in \widetilde{\mathcal{G}}}$ is $L^{1}(P)$-bounded and claim i) follows.

\emph{ii) }Having proved i) $\{U\left( X^{+}\right) \}_{X\in \widetilde{\mathcal{G}}}$ is $L^{1}(P)$-bounded 
by (\ref{eq: Fenchel Uplus}), and hence
by the assumed lower bound on $\{I_{U}(X)\}_{X\in \widetilde{\mathcal{G}}}$ the set 
$\{U(-X^{-})\}_{X\in \widetilde{\mathcal{G}}}$ is also $L^{1}(P)$-bounded .

\emph{iii) }Item ii) implies norm-boundedness of $\{X^{-}\}_{X\in \widetilde{\mathcal{G}}}$ 
by equivalence of gauge norms \citep[Proposition~2]{caruso.01}. 
Item iii) now follows by contradiction because $\widetilde{\mathcal{G}}$ was arbitrary.
\end{proof}

\begin{corollary}
\label{cor: XY L(P)-bounded}Suppose $\{U(X)^{-}\}_{X\in \widetilde{\mathcal{G}}}$ 
is $L^{1}(P)$-bounded. For any $Y\in L^{\hat{V}}$ such that 
$\{E[XY]\}_{X\in \widetilde{\mathcal{G}}}$ is bounded above 
$\{\left\vert X\right\vert Y\}_{X\in \widetilde{\mathcal{G}}}$ is also $L^{1}(P)$-bounded.
\end{corollary}

\begin{proof}
$L^{1}(P)$-boundedness of $\{X^{-}Y\}_{X\in \widetilde{\mathcal{G}}}$
follows from the Fenchel inequality (\ref{eq: Fenchel Uminus}) where we take 
$\lambda _{2}$ such that $\lambda _{2}Y\in \mathrm{dom}\,I_{\hat{V}}$. $L^{1}(P)$-boundedness of 
$\{X^{+}Y\}_{X\in \widetilde{\mathcal{G}}}$ now
follows from the assumed upper bound on $\{E[XY]\}_{X\in \widetilde{\mathcal{G}}}$.
\end{proof}

As the final ingredient we must ensure uniform integrability of $\{U(A_{n}^{+})\}$ 
for a maximizing sequence $\{A_{n}\}=B+\{X_{n}\}$. Define
indirect utility $\overline{u}:\mathbb{R}_{+}\rightarrow \mathbb{R}$ by
maximizing $u(B)$ over all random endowments $B$ whose $L^{\hat{U}}$ norm is
bounded above by $x$, 
\begin{equation}
\overline{u}(x)=\sup_{\Vert B\Vert _{\hat{U}}\leq x}u(B)=\sup
\{I_{U}(X+Z)\mid X\in \mathcal{C},\Vert Z\Vert _{\hat{U}}\leq x\}.
\label{eq: ubar def}
\end{equation}
To obtain uniform integrability of positive parts of utility we will require
\begin{equation}
\lim_{x\rightarrow \infty }\overline{u}(x)/x=0.  \label{eq: Inada ubar}
\end{equation}
The construction involving $\overline{u}$ also appears to be new.

Note that for $L^{\hat{U}}\sim L^{\infty }$ one has $\overline{u}(x)=
u(x)$ for $x\geq 0$ and therefore condition (\ref{eq: Inada ubar}) exactly
coincides with the minimal condition in \citet[Note~1]{kramkov.schachermayer.03}. 
The significance of the condition (\ref{eq: Inada ubar}) is captured by
the following statement.

\begin{lemma}
\label{lem: UI U(mathcal Z)}Condition $\lim_{x\rightarrow \infty }\overline{u}(x)/x=0$ implies that for
\begin{equation*}
\mathcal{Z}(k_{1},k_{2})=\{X+Z\mid X\in \mathcal{C},Z\in L^{\hat{U}},\Vert
X^{-}\Vert _{\hat{U}}\leq k_{1},\Vert Z\Vert _{\hat{U}}\leq k_{2}\},
\end{equation*}
the set $\{U(\mathcal{Z}(k_{1},k_{2})^{+})\}$ is uniformly integrable for every $k_{1},k_{2}>0.$
\end{lemma}

\begin{proof}
It follows from the Eberlein-\v{S}mulian and Dunford-Pettis theorems 
\citep[Theorems~4.7.10 and 4.7.18]{bogachev.07}
that uniform integrability (UI) of a set is equivalent to UI of sequences in
the set (see also \citealp[pages~45 and 50]{diestel.91}). Now we can proceed as in 
\citet{kramkov.schachermayer.03}
but with the new notion $\overline{u}$ in place of $u$ which allows us to
handle the general case where the unit ball of $L^{\hat{U}}$ does not have
an upper bound. We also replace polarity arguments of the original proof
(unavailable here) with simpler set inclusions.

Arguing by contradiction assume that for $X_{i}\in \mathcal{C},\Vert
X_{i}^{-}\Vert _{\hat{U}}\leq k_{1}$ and $\Vert Z_{i}\Vert
_{\hat{U}}\leq k_{2}$ the sequence $\{U((X_{i}+Z_{i})^{+})\}$ is not
uniformly integrable. Then there are disjoint sets $D_{i}\in \mathcal{F}_{T}$
and a constant $\alpha >0$ such that 
\begin{equation*}
E[U((X_{i}^{+}+Z_{i}^{+})1_{D_{i}})]\geq
E[U((X_{i}+Z_{i})^{+}1_{D_{i}})]\geq \alpha .
\end{equation*}
Note that for $X\in \mathcal{C},\Vert X^{-}\Vert _{\hat{U}}\leq
k_{1}$ one has
\begin{equation*}
X_{i}^{+}=X_{i}+X_{i}^{-}\in \mathcal{Z}(k_{1},k_{1}),
\end{equation*}
which implies $$\sum_{i=1}^{n}(X_{i}^{+}+Z_{i}^{+})1_{D_{i}}\leq
\sum_{i=1}^{n}(X_{i}^{+}+Z_{i}^{+})\in \mathcal{Z}(nk_{1},n(k_{1}+k_{2}))$$
and consequently
\begin{equation*}
n\alpha \leq
\sum_{i=1}^{n}I_{U}(X_{i}^{+}+Z_{i}^{+})1_{D_{i}}
=I_{U}\left(\sum_{i=1}^{n}(X_{i}^{+}+Z_{i}^{+})1_{D_{i}}\right)\leq \overline{u}(n(k_{1}+k_{2})).
\end{equation*}
From here $\alpha /(k_{1}+k_{2})\leq \overline{u}(n(k_{1}+k_{2}))/\left(n(k_{1}+k_{2})\right)$, 
and for $n\rightarrow \infty $ the right-hand side
converges to $0$ by hypothesis which gives the desired contradiction.
\end{proof}

We are now in a position to prove the existence of an optimal terminal wealth with
the desired approximation property.

\begin{proposition}\label{prop: U(X_n) converges to U(hatX) in L1(P)}
Assume $B\in L^{\hat{U}}$
and $(B+\mathcal{C})\cap \mathrm{core\,}\mathcal{D\neq \emptyset }$. Assume
further there is no arbitrage over $\mathcal{C}^{\circledast \circledast }$; 
$\lim_{x\rightarrow \infty }\overline{u}(x)/x=0$; and the dual minimizer $\hat{Y}$ in 
(\ref{eq: duality AcapD}) satisfies $\lambda \hat{Y}\in \mathrm{dom}\,I_{V}$ for
some $\lambda >1$ (this is automatic when $M^{\hat{V}}=L^{\hat{V}}$). Then
there is a sequence $\{X_{n}\}\in \mathcal{C}$ with $I_{U}(B+X_{n})\nearrow
u(B)<\infty $ and a random variable $\hat{X}$ such that $X_{n}\overset{P\text{-a.s.}}{\rightarrow }\hat{X}$, 
\begin{equation}
U(B+\hat{X})-(B+\hat{X})\hat{Y}=V(\hat{Y}),
\label{eq: Fenchel equality Xhat Yhat}
\end{equation}
and $$U(B+X_{n})\overset{L^{1}(P)}{\rightarrow }U(B+\hat{X}).$$ Moreover, the sequence $\{X_n\}$ can be chosen such that 
$U(B+X_n)\geq R$ with $0\geq R\in L^1(P)$.
\end{proposition}

\begin{proof}
\emph{Step 1)} We will first exhibit a random variable $\tilde{Y}>0$ $P$-a.s. 
such that $\lambda \tilde{Y}\in \mathrm{dom}\,I_{V}$ for two distinct
values of $\lambda $ and $\sup_{X\in \mathcal{A}\cap \mathcal{D}}\{E[X\tilde{Y}]\}<\infty$. 
We distinguish two mutually exclusive cases.

\emph{a)} When $U$ is bounded from above then $V(y)$ is bounded from above
for $y$ near zero. By the Kreps-Yan theorem \citep[Proposition~3.5]{gao.xanthos.17}
no arbitrage over $\mathcal{C}^{\circledast \circledast }$ implies
existence of $\tilde{Y}\in \mathcal{C}^{\circledast },\tilde{Y}>0$ $P$-a.s.
Because $\mathcal{C}^{\circledast }$ is a cone, without loss of generality we
may assume $\tilde{Y}\in \mathrm{dom}\,I_{\hat{V}}$. Recalling that $\hat{V}(y)=V(\vert y\vert\vee a)$ 
while $V(\vert y\vert\wedge a)$ is bounded we conclude $\lambda \tilde{Y}\in \mathrm{dom}\,I_{V}$ 
for all $0<\lambda \leq 1$. By construction $\sup_{X\in \mathcal{C}}\{E[X\tilde{Y}]\}\leq 0$ which implies $\sup_{X\in 
\mathcal{A}}\{E[X\tilde{Y}]\}<\infty $.

\emph{b)} By condition (\ref{eq: Inada ubar}) $u(B)<\infty $. When $U$ is
unbounded from above then $V(0)=\infty $ and therefore necessarily the dual
optimizer in (\ref{eq: duality AcapD}) satisfies $\hat{Y}>0$ $P$-a.s. as
well as $\sup_{X\in \mathcal{A}\cap \mathcal{D}}\{E[X\hat{Y}]\}<\infty $ and 
$\hat{Y}\in \mathrm{dom}\,I_{V}.$ In this case we let $\tilde{Y}=\hat{Y}$.

\emph{Step 2)} By definition of supremum there is a sequence $\{A_{n}\}$ in 
$\mathcal{A}$ with $\{I_{U}(A_{n})\}$ bounded below and $I_{U}(A_{n})\nearrow
u(B)$. The random variable $\tilde{Y}$ from step 1) and the sets $\mathcal{G}=\mathcal{A}\cap \mathcal{D}$ 
and $\widetilde{\mathcal{G}}=\mathrm{conv}\{A_{n}\}$ therefore satisfy the hypotheses of 
Lemma~\ref{lem: L^1(P) boundedness}. We thus conclude that 
$\{U(A_{n}^{+})\}$, $\{U\left(-A_{n}^{-}\right) \}$, $\{A_{n}^{+}\tilde{Y}\}$, $\{A_{n}^{-}\tilde{Y}\}$
are $L^{1}(P)$-bounded.

\emph{Step 3)}\ Construct $\hat{A}\tilde{Y}$ as the pointwise limit of tail
convex combinations of $A_{n}\tilde{Y}$. By abuse of notation denote these
convex combinations again $A_{n}\tilde{Y}$. By construction $A_{n}\rightarrow \hat{A}$ $P$-a.s. 
Note that the new sequence $\{A_{n}\}$
satisfies the same hypotheses as the old one: all elements are in $\mathcal{A}$ and $I_{U}(A_{n})$ 
is bounded from below and converges to $u(B)$. Let $X_{n}=A_{n}-B$ and $\hat{X}=\hat{A}-B$.

\emph{Step 4)} From here onwards we pass to a subsequence such that 
$I_{U}(A_{n}^{+})$, $I_{U}(-A_{n}^{-})$, $E[A_{n}^{+}\hat{Y}]$, and
$E[A_{n}^{-}\hat{Y}]$ all have a finite limit.

\emph{Step 5)} By assumption there is $\lambda _{0}>1$ such that for all $\lambda \in [ 1,\lambda _{0}]$ 
we have $\lambda \hat{Y}\in \mathrm{dom}\,I_{V}$. Fatou lemma yields
\begin{align*}
\lim_{n\rightarrow \infty }I_{U}(A_{n})-\lambda \lim_{n\rightarrow \infty
}E[A_{n}\hat{Y}] &=\lim_{n\rightarrow \infty }\{I_{U}(A_{n})-\lambda E[A_{n}\hat{Y}]\} \\
&\leq I_{U}(\hat{A})-\lambda E[\hat{A}\hat{Y}]\leq I_{V}(\lambda \hat{Y}),
\end{align*}
which means
\begin{equation}\label{eq: Fatou inequality lambda}
u(B)-\lambda \lim_{n\rightarrow \infty }E[A_{n}\hat{Y}]
\leq I_{U}(\hat{A})-\lambda E[\hat{A}\hat{Y}]\leq I_{V}(\lambda \hat{Y}).
\end{equation}

\emph{Step 6)} By Theorem~\ref{thm: strong duality} 
\begin{equation}
u(B)=I_{V}(\hat{Y})+\sup_{A\in \mathcal{A}\cap \mathcal{D}}E[A\hat{Y}].
\label{propXhat1}
\end{equation}
Substitute this into (\ref{eq: Fatou inequality lambda}) with $\lambda =1$
to obtain
\begin{equation*}
I_{V}(\hat{Y})+\sup_{A\in \mathcal{A}\cap \mathcal{D}}E[A\hat{Y}]-\lim_{n\rightarrow \infty }E[A_{n}\hat{Y}]
\leq I_{V}(\hat{Y}).
\end{equation*}
This implies $\sup_{A\in \mathcal{A}\cap \mathcal{D}}E[A\hat{Y}]-\lim_{n}E[A_{n}\hat{Y}]\leq 0$ 
but as $A_{n}\in \mathcal{A}\cap\mathcal{D}$, this is only possible if 
\begin{equation}
\lim_{n\rightarrow \infty }E[A_{n}\hat{Y}]=\sup_{A\in \mathcal{A}\cap 
\mathcal{D}}E[A\hat{Y}].  \label{propXhat2}
\end{equation}
Therefore for $\lambda =1$ the inequalities in (\ref{eq: Fatou inequality
lambda}) are actually equalities 
\begin{equation}
u(B)-\lim E[A_{n}\hat{Y}]=I_{U}(\hat{A})-E[\hat{A}\hat{Y}]=I_{V}(\hat{Y}).  \label{propXhat3}
\end{equation}
Equality (\ref{propXhat3}) implies that the Fenchel inequality $U(\hat{A})-\hat{A}\hat{Y}\leq V(\hat{Y})$ 
is in fact a $P$-a.s. equality which proves (\ref{eq: Fenchel equality Xhat Yhat}).

\emph{Step 7)} Subtract (\ref{propXhat3}) from (\ref{eq: Fatou inequality lambda}) to obtain 
\begin{equation*}
(1-\lambda )\lim_{n\rightarrow \infty }E[A_{n}\hat{Y}]
\leq (1-\lambda )E[\hat{A}\hat{Y}]\leq I_{V}(\lambda \hat{Y})-I_{V}(\hat{Y}).
\end{equation*}
Taking $\lambda >1$ we have
\begin{equation}
\lim_{n\rightarrow \infty }E[A_{n}\hat{Y}]\geq E[\hat{A}\hat{Y}].
\label{propXhat4}
\end{equation}
On combining (\ref{propXhat4}) with (\ref{propXhat1}) and (\ref{propXhat2})
we finally conclude
\begin{equation}
u(B)\geq I_{U}(\hat{A})=I_{V}(\hat{Y})+E[\hat{A}\hat{Y}].
\label{propXhat5}
\end{equation}

\emph{Step 8)} Observe that a sequence $Z_{n}\overset{P\text{-a.s.}}{\rightarrow }Z$ is uniformly integrable if and only if 
$$E[\left\vert Z_{n}\right\vert]\rightarrow E[\left\vert Z\right\vert ] \iff 
E[\left\vert Z_{n}-Z\right\vert ]\rightarrow 0,$$ see Scheff\'{e} lemma \citep[Theorem~2.8.9]{bogachev.07}
and Lebesgue-Vitali convergence theorem \citep[Theorem~4.5.4]{bogachev.07}. 
Fatou lemma yields 
\begin{equation}
\lim_{n}I_{U}(-A_{n}^{-})\leq I_{U}(-\hat{A}^{-}).
\label{eq: Fatou I_U minus}
\end{equation}
In order to obtain $L^{1}(P)$-convergence of $\{I_{U}(A_{n})\}$ in view of (\ref{propXhat5}) 
and (\ref{eq: Fatou I_U minus}) it suffices to prove $I_{U}(A_{n}^{+})\rightarrow I_{U}(\hat{A}^{+})$ 
or equivalently that the sequence $\{U(A_{n}^{+})\}$ is uniformly integrable.

\emph{Step 9)} By step 2) $\sup_{n}\Vert (B+X_{n})^{-}\Vert _{\hat{U}}<\infty$. Let $k_{1}=\left\Vert B\right\Vert _{\hat{U}}$. 
We have $X_{i}^{-}\leq (B+X_{i})^{-}+B^{+}$ and therefore 
$\Vert X_{i}^{-}\Vert _{\hat{U}}\leq \sup_{n}\Vert (B+X_{n})^{-}\Vert _{\hat{U}}+\Vert B\Vert _{\hat{U}}=:k_{2}<\infty$. 
We conclude that $A_{n}=B+X_{n}\in \mathcal{Z}(k_{1},k_{2})$ and the sequence $\{U(A_{n}^{+})\}$ is uniformly integrable
by Lemma~\ref{lem: UI U(mathcal Z)}. By step 8) $U(A_{n})\rightarrow U(\hat{A})$ in $L^{1}(P)$.

\emph{Step 10)} This means the non-positive sequence $\{U(-A_{n}^{-})\}$ 
is Cauchy in $L^{1}(P)$ and we can find a subsequence, here denoted by $\tilde{A}_{n}$, and a random variable 
$R\in L^{1}(P)$ such that $0\geq U(-\tilde{A}_{n}^{-})\geq R$.
\end{proof}

\begin{remark}
Corollary~3.10 in \citet{delbaen.owari.16} shows, under the assumption that $\hat{V}$ satisfies the $\Delta_2$--condition, 
that every $L^{\hat{U}}$ norm-bounded sequence admits a pointwise-convergent sequence of forward convex combinations 
whose $\hat{U}$ is dominated by an integrable random variable. Here the dominated convergence of forward convex 
combinations is shown to exist for the negative parts of the sequence of terminal wealths, $\{A_n^{-}\}$, without 
necessarily assuming the $\Delta_2$--condition on $\hat{V}$. Our starting sequence, however, is maximizing and therefore not arbitrary.
\end{remark}

\section{Duality over separating measures}\label{sect: duality over separating}

We have argued in the
introductory Section~\ref{sect: intro optimal strategy} that for the
construction of the optimal trading strategy it is important to know there
exists a full market completion whose utility is arbitrarily close to $u(B)$,
\begin{equation}\label{eq: duality over separating meas}
u(B)=\sup_{X\in \mathcal{C}}I_{U}(B+X)=\inf_{Y\in \mathcal{C}^{\circledast}}\{I_{V}(Y)+E[YB]\}.  
\end{equation}

The case $M^{\hat{U}}=L^{\hat{U}}$ is immediately very nice in this respect:
one automatically has $\star =\circledast $ so the norm duality (\ref{eq: norm duality}) 
yields the desired result (\ref{eq: duality over separating meas}). This
covers the case \textbf{L }where the utility function is asymptotically
linear near $-\infty $ and $L^{\hat{U}}\sim L^{1}$.

The remaining case is \textbf{SL} with $M^{\hat U}\subsetneq L^{\hat U}$.  The norm 
duality (\ref{eq: norm duality}) implies that utility cannot 
increase by going from $B+\mathcal{C}$ to its norm-closure $B+\mathcal{C}^{\star \star }$
while the economic duality (\ref{eq: dual corner 2}) implies that utility does
not increase by going from $(B+\mathcal{C})\cap \mathcal{D}$ to $\mathrm{cl}^{\circledast }((B+\mathcal{C})\cap \mathcal{D})$. 
However, these facts do not in themselves prevent a utility gap between 
$B+\mathcal{C}$ and $B+\mathcal{C}^{\circledast \circledast }$. 
Appendix~\ref{sect: utility increases} provides a counterexample illustrating
that with $M^{\hat{U}}\subsetneq L^{\hat{U}}$ one can generically expect to find situations where
\begin{equation}\label{eq: utility increases}
\sup_{X\in \mathcal{C}}I_U(B+X) < \sup_{X\in \mathcal{C}^{\circledast\circledast}}I_U(B+X).
\end{equation}
Obviously, if the gap \eqref{eq: utility increases} emerges then by the Fenchel inequality \eqref{eq: duality over separating meas} cannot hold.

This observation highlights the importance of the classical
`small market' fundamental theorem of asset pricing (FTAP) which asserts 
that in the absence of arbitrage over $\mathcal{C}^{\star\star}$ in the case \textbf{INF} 
one necessarily obtains $\mathcal{C}=\mathcal{C}^{\circledast\circledast}$. 
Our counterexample also shows that in a `large financial market' the link 
between absence of arbitrage over $\mathcal{C}^{\circledast \circledast }$ and 
the equality $\mathcal{C}=\mathcal{C}^{\circledast \circledast }$ is broken, 
and the case \textbf{INF} is no exception.

Having made the necessary preparations, it turns out that the following weaker 
alternative of \eqref{eq: duality over separating meas} is already sufficient for our purposes.  

\begin{proposition}\label{prop: y + Y}
Assume $B\in L^{\hat{U}}$,
$(B+\mathcal{C})\cap \mathcal{D\neq \emptyset }$, and 
\begin{equation}\label{eq: utilities equal}
u(B)=\sup_{X\in \mathcal{C}}I_U(B+X) = \sup_{X\in \mathcal{C}^{\ast\ast}}I_U(B+X).
\end{equation}
 Then in the case \textbf{SL} ($L^{\infty }(P)
\hookrightarrow L^{\hat{U}}(P)
\hookrightdoublearrow L^{1}(P)$),  there is a
sequence $\{Z_{n}\}\in M^{\hat{V}}$ with $\left\Vert Z_{n}\right\Vert _{\hat{V}}\rightarrow 0$, 
and a sequence of $\{Y_{n}\}\in \mathcal{C}^{\ast }$ such that 
$$\lim_{n\rightarrow \infty }I_{V}\left(Y_{n}+Z_{n}\right)+E[(Y_n+Z_n) B] =u(B).$$
\end{proposition}

\begin{proof}
By Proposition~\ref{prop: I_U is usc} $I_{U}$ is $\ast $-u.s.c. Because $I_{U}$ is 
finite-valued at $0$ it is proper by Proposition~\ref{prop: usc h is proper}, 
and therefore $\ast $-closed by Definition~\ref{def: closure}.
Likewise $\delta _{\mathcal{C}^{\ast \ast }}$ is a $\ast $-closed function because $\mathcal{C}^{\ast \ast }$ 
is a closed (convex) set in the duality $(L^{\hat{U}},M^{\hat{V}})$. 
Taking $f(X) = I_{U}(B+X)$ and $g= -\delta _{\mathcal{C}^{\ast \ast }}$
we have $f+g$ is $\ast $-u.s.c. by Proposition~\ref{lem: sum of usc is usc}.
The sum is also proper and therefore closed because $ \mathrm{dom}\,f\cap \mathrm{dom}\,g\neq\emptyset$. We have $f^\ast(Y)=-I_V(Y)-E[YB]$ and $g^\ast = -\delta_{-\mathcal{C}^\ast}$. 
Because $f^\ast(1)$ and $g^\ast(0)$ are finite, 
 $f^\ast,g^\ast$ are proper and by Lemma~\ref{lem: (f+g)*}
\begin{equation}\label{eq: lsc square}
u(B)=\sup_{X\in L^{\hat{U}}}\{I_{U}(B+X)-\delta _{\mathcal{C}^{\ast \ast }}(X)\}
=\mathrm{lsc}\,(-f^\ast\conv \delta _{-\mathcal{C}^{\ast }})(0).
\end{equation}
Due to $L^{\hat{U}}=(M^{\hat{V}})^{\star }$ we may evaluate the lower
semicontinuous hull in the norm topology on $M^{\hat{V}}$, see Theorem~\ref{thm: lsc hull invariance}. 
Therefore there exists a sequence $Z_{n}$ in $M^{\hat{V}}$ norm-convergent to $0$ such that 
$\lim_{n\rightarrow \infty}(-f^\ast\conv \delta _{-\mathcal{C}^{\ast }})(Z_{n})=u(B)$ which completes the proof on
recalling the formula for the infimal convolution, see Definition~\ref{def: inf convolution}.
\end{proof}

\begin{remark}
In the case $M^{\hat{V}}=L^{\hat{V}}$ one has $\circledast=\ast$ hence the 
assumption \eqref{eq: utilities equal} is absolutely necessary to prevent 
the utility gap in \eqref{eq: utility increases}. In contrast, with $M^{\hat{V}}\subsetneq L^{\hat{V}}$
condition \eqref{eq: utilities equal}  is no longer economically innocuous because
there are complete market examples where $\mathcal{C}^{\ast \ast }=L^{\hat{U}}$ while 
$\mathcal{C}^{\circledast \circledast }$ is arbitrage-free. Nonetheless, assumption 
\eqref{eq: utilities equal} gives, by some margin, the best 
result available to date. 

At present the only works in the literature that allow 
$M^{\hat{V}}\subsetneq L^{\hat{V}}$ are \citet{biagini.frittelli.05} and 
\citet{biagini.cerny.11} who require $\mathcal{C}^{\circledast }=\mathcal{C}^{\ast }$ 
which forces $\mathcal{C}^{\circledast \circledast }=\mathcal{C}^{\ast \ast }$ and so implies 
\eqref{eq: utilities equal}. 
\citet{biagini.frittelli.08} assume $M^{\hat{V}}=L^{\hat{V}}$ and 
\citet{biagini.frittelli.07}, \citet{schachermayer.01,schachermayer.03},
and \citet{owen.zitkovic.09} require reasonable asymptotic elasticity at $-\infty$
\citep[Definition~1.4]{schachermayer.01} which is stronger than $M^{\hat{V}}=L^{\hat{V}}$ 
\citep[Proposition~4.1(iii)]{schachermayer.01}.
\end{remark}

For completeness we now prove the full duality over separating measures 
(\ref{eq: duality over separating meas}) which requires stronger assumptions.

\begin{theorem}
\label{thm: duality over MM}Assume either i) $M^{\hat{U}}=L^{\hat{U}}$; or
ii) $\mathcal{C}=\mathcal{C}^{\ast \ast }$; $\lim_{x\rightarrow \infty}\overline{u}(x)/x=0$; 
there is $0<\bar{Y}\in \mathcal{C}^{\circledast }$
(no arbitrage over $\mathcal{C}^{\circledast \circledast }$); and, only in
the case \textbf{F-SL}, there is $\tilde{Y}\in \mathcal{C}^{\circledast }$
such that $\lambda \tilde{Y}\in \mathrm{dom}\,I_{V}$ for two distinct values
of $\lambda \geq 0$. Then the duality over separating measures (\ref{eq: duality
over separating meas}) holds for all $B\in L^{\hat{U}}$.
\end{theorem}

\begin{proof}
\emph{Step 1)}\ For $M^{\hat{U}}=L^{\hat{U}}$ the claim follows from 
Theorem~\ref{thm: strong duality}. This covers case \textbf{L}. It remains to prove
the case \textbf{SL} under the assumption ii). Recall $u:L^{\hat{U}}\rightarrow\mathbb{R}\cup \{-\infty \}$ 
is the maximal expected utility as a function of the random endowment $Z\in L^{\hat{U}}$,
\begin{equation*}
u(Z)=\sup_{X\in \mathcal{C}}\{I_{U}(X+Z)\}=\left( I_{U}\conv -\delta _{-\mathcal{C}}\right) (Z),
\end{equation*}
where $\conv $ denotes the supremal convolution (Definition~\ref{def: inf convolution}). 

Because both $I_{U}$ and $-\delta _{-\mathcal{C}}$ are proper
(Definition~\ref{def: effective domain}), by Lemma~\ref{lem: (f+g)*} 
$$u^{\ast }(Y)=I_{U}^{\ast }(Y)-\delta _{\mathcal{C}}^{\ast }(Y)=-I_{V}(Y)-\delta _{\mathcal{C}^{\ast }}(Y)$$ 
and by the definition of the conjugate function
\begin{equation}\label{eq: h**}
u^{\ast \ast }(Z)=\inf_{Y\in M^{\hat{V}}}\{E[YZ]+I_{V}(Y)+\delta _{\mathcal{C}^{\ast }}(Y)\}
=\inf_{Y\in \mathcal{C}^{\ast }}\{E[YZ]+I_{V}(Y)\}.
\end{equation}

\emph{Step 2) }By virtue of (\ref{eq: h**}) the proof will be complete if we
can show $u(B)=u^{\ast \ast }(B)$. By Proposition~\ref{prop: usc h is proper}
and Theorem~\ref{thm: Fenchel-Moreau} this is equivalent to demonstrating
that $u$ is $\ast $-u.s.c. at $B$. This line of reasoning is the essence of
the conjugate duality construction proposed in \citet{rockafellar.74}. 
We will show a stronger property, namely that $u$ is $\ast $-u.s.c.
globally. By Proposition~\ref{prop: gao.xanthos.18} $u$ is $\ast $-u.s.c. if
and only if for arbitrary norm-bounded sequence $\{Z_{n}\}\in L^{\hat{U}}$
such that $Z_{n}\overset{P\text{-a.s.}}{\rightarrow }Z\in L^{\hat{U}}$ one
has $\underset{n\rightarrow \infty }{\lim \sup }\,u(Z_{n})\leq u(Z)$.

\emph{Step 3)\ }If $\lim \sup_{n\rightarrow \infty }\,u(Z_{n})=-\infty $
there is nothing to prove. In the remaining case $\lim \sup_{n\rightarrow
\infty }\,u(Z_{n})=:\tilde{u}>-\infty $. By the definition of supremum there is
a subsequence (still denoted $Z_{n}\in L^{\hat{U}}$) and a corresponding
sequence of $X_{n}\in \mathcal{C}$ such that $I_{U}(X_{n}+Z_{n})$ is bounded
below and 
\begin{equation*}
I_{U}(X_{n}+Z_{n})\nearrow \tilde{u}.
\end{equation*}
Denote by $\widetilde{\mathcal{G}}$ the convex hull of $\{X_{n}+Z_{n}\}$. By
convexity of upper level sets $I_{U}$ is bounded below on $\widetilde{\mathcal{G}}$.

\emph{Step 4)} We claim that $k\mathcal{B}+\mathcal{C}$ is norm-coercive in
losses (see Definition~\ref{def: norm-coercive in losses}) for arbitrary $k>0 $, where $\mathcal{B}$ 
is the unit ball in $L^{\hat{U}}$. For $L^{\hat{U}}\sim L^{\infty }$ 
this is true trivially. In the remaining case \textbf{F-SL} the set 
$\mathcal{G}=k\mathcal{B}+\mathcal{C}$ and separating density $\tilde{Y}$ satisfy the
assumptions of Lemma~\ref{lem: L^1(P) boundedness} and
the claim follows. As a result $\widetilde{\mathcal{G}}$ is norm-bounded in
losses and in view of the norm-boundedness of $\{Z_{n}\}$ the set $\{U(\widetilde{\mathcal{G}}^{+})\}$ 
is uniformly integrable by Lemma~\ref{lem: UI U(mathcal Z)}. By Lemma~\ref{lem: L^1(P) boundedness} 
$\{\left\vert X_{n}+Z_{n}\right\vert \tilde{Y}\}$ is $L^{1}(P)$-bounded, while 
$\{\left\vert X_{n}\right\vert \tilde{Y}\}$, too, is $L^{1}(P)$-bounded by the H\"{o}lder inequality. 
The same holds for $\{\left\vert X_{n}+Z_{n}\right\vert \bar{Y}\}$, $\{\left\vert X_{n}\right\vert \bar{Y}\}$
by Corollary \ref{cor: XY L(P)-bounded}.

\emph{Step 5)} Letting $d\bar{Q}=\bar{Y}dP$, Koml\'{o}s theorem yields a
sequence of forward convex combinations of $\{X_{n}\}$ (denoted $\{\hat{X}_{n}\}$) such that 
$\hat{X}_{n}$ converges $P$-a.s. to some limit $\hat{X}$.
We will apply the same convex combinations to $\{Z_{n}\}$ and denote the
resulting sequence by $\{\hat{Z}_{n}\}$. By construction $\hat{X}_{n}+\hat{Z}_{n}\in \widetilde{\mathcal{G}}$, 
therefore $U((\hat{X}_{n}+\hat{Z}_{n})^{+})$ is uniformly integrable by step 4). By concavity of $I_{U}$ the
utility of convex combinations dominates the utility of the original sequence, 
\begin{equation*}
\tilde{u}\leq \hat{u}=\lim \sup I_{U}(\hat{X}_{n}+\hat{Z}_{n}),
\end{equation*}
and by passing to a further subsequence we may assume $\hat{u}=\lim I_{U}(\hat{X}_{n}+\hat{Z}_{n})$.

\emph{Step 6) }UI of $\{U((\hat{X}_{n}+\hat{Z}_{n})^{+})\}$ and Fatou lemma yield 
\begin{align}
\hat{u} &=\lim_{n\rightarrow \infty }I_{U}(\hat{X}_{n}+\hat{Z}_{n})
=\lim_{k\rightarrow \infty }\underset{n\rightarrow \infty }{\lim \sup }I_{U}((\hat{X}_{n}+\hat{Z}_{n})\wedge k)  \notag \\
&\leq \lim_{k\rightarrow \infty }I_{U}((\hat{X}+Z)\wedge k).\label{eq: Fatou k}
\end{align}
For any $k>0$ the sequence $\{\hat{X}_{n}\wedge (k-Z_{n})\}\in \mathcal{C}$
is norm-bounded and $P$-a.s. convergent to $\hat{X}\wedge (k-Z)$. By 
\citet[Theorem~2.1]{gao.14}
we conclude that $\hat{X}\wedge (k-Z)\in \mathcal{C}^{\ast \ast }$. By
assumption $\mathcal{C}=\mathcal{C}^{\ast \ast },$ therefore 
\begin{equation*}
I_{U}((\hat{X}+Z)\wedge k))=I_{U}(Z+\hat{X}\wedge (k-Z))\leq \sup_{X\in \mathcal{C}}I_{U}(Z+X)= u(Z),
\end{equation*}
and from (\ref{eq: Fatou k}) we conclude $\hat{u}\leq u(Z)$.

\emph{Step 7) }Combining the steps 1)-6) we have shown in the case \textbf{SL}
under the assumption ii) 
\begin{equation*}
\sup_{X\in \mathcal{C}}\{I_{U}(X+Z)\}=\inf_{Y\in \mathcal{C}^{\ast
}}\{I_{V}(Y)+E[YZ]\}\text{ for all }Z\in L^{\hat{U}}\text{.}
\end{equation*}
By the Fenchel inequality $\sup_{X\in \mathcal{C}}\{I_{U}(X+Z)\}\leq \inf_{Y\in 
\mathcal{C}^{\circledast }}\{I_{V}(Y)+E[YZ]\}$ while 
\begin{equation*}
\inf_{Y\in \mathcal{C}^{\circledast }}\{I_{V}(Y)+E[YZ]\}\leq \inf_{Y\in 
\mathcal{C}^{\ast }}\{I_{V}(Y)+E[YZ]\},
\end{equation*}
due to the inclusion $\mathcal{C}^{\ast }=\mathcal{C}^{\circledast }\cap M^{\hat{V}}\subseteq \mathcal{C}^{\circledast }$. 
This proves (\ref{eq: duality over separating meas}).
\end{proof}

\section{Optimal admissible strategy}\label{sect: optimal strategy}

Our soon-to-be-found ability to deal with utility functions that are not
strictly monotone prompts a slight modification of the definition of
admissibility, compared to \citet[Definition~1.1]{biagini.cerny.11}. 
In this paper we require a tight approximation of the wealth process below
the bliss point of the utility function but only a loose one above the bliss
point.

We also need to amend the definition of convergence at intermediate times to
allow for effective completions as dual optimizers. The limiting process has
to be defined not as a pointwise limit of $H^{(n)}\cdot S$ at fixed times
but rather as its right-continuous regularization\footnote{For regularization of submartingales see, 
for example, \citet[Theorem~II.2.5 and Proposition~II.2.6]{revuz.yor.99}} 
$\mathrm{rqlim}_{n\rightarrow \infty }H^{(n)}\cdot S$ defined by
\begin{equation}\label{eq: rqlim}
(\mathrm{rqlim}_{n\rightarrow \infty }H^{(n)}\cdot S)_{t}=\lim_{q\searrow
t,q\in \mathbb{Q}}(\lim_{n\rightarrow \infty }H^{(n)}\cdot S_{q}),
\end{equation}
in line with the supermartingale compactness result in \citet[Theorem~D]{delbaen.schachermayer.98}.

\begin{definition}
\label{def: admissible}A strategy $H\in L(S)$ is admissible 
($H\in \scr{A}$) if there is a sequence $H^{(n)}$ in $\scr{T}$ such that
\begin{enumerate}
\item $U\left( B+H^{(n)}\cdot S_{T}\right) \rightarrow U\left( B+H\cdot S_{T}\right) $ in $L^{1}(P)$;

\item on the set $U(B+H\cdot S_{T})<U(\infty )$ the approximating tame
wealth $H^{(n)}\cdot S$ converges to the admissible wealth $H\cdot S$ in the
sense of right-continuous regularization (\ref{eq: rqlim}) 
\begin{equation*}
H\cdot S_{t}=(\mathrm{rqlim}_{n\rightarrow \infty }H^{(n)}\cdot S)_{t}\qquad\text{for all\ }t\in [ 0,T].
\end{equation*}
\end{enumerate}
\end{definition}

Recall the definition of the set of supermartingale measures $\mathcal{S}$ in equation 
(\ref{eq: supermart measures}). We begin by observing that the wealth process of every tame
strategy is a supermartingale under each $Q\in \mathcal{M}\cap P_{\hat{V}}.$

\begin{proposition}
\label{prop: tame strategies are super}For $Q\in \mathcal{M}\cap P_{\hat{V}}$
and $H\in \scr{T}$ the wealth process $H\cdot S$ is a $Q$-supermartingale. 
In other words, $\mathcal{M}\cap P_{\hat{V}}\subseteq \mathcal{S}$.
\end{proposition}

\begin{proof}
\emph{Step 1) }$S$ is a $Q$-$\sigma $-martingale therefore $H\cdot S$ can be
written as an integral with respect to a $Q$-martingale \citep[Proposition~2]{emery.80}. 
$H\in \scr{T}$ means that for $W_{t}=\inf_{\tau \in [ 0,t]}\{H\cdot S_{\tau
}\}\wedge 0$ one has $W_{T}\in L^{\hat{U}}$. This implies $W_{T}\in L^{1}(Q)$
and there is a $P$-martingale $Z^{Q}$ such that $Z_{T}^{Q}=W_{T}$. The
minimal process $W$ is decreasing and therefore $W\geq Z^{Q}$.

\emph{Step 2)} Because $H\cdot S$ is bounded below by the $Q$-martingale $Z^{Q}$, 
it follows from the Ansel-Stricker lemma \cite[Corollaire~3.5]{ansel.stricker.94} 
that $H\cdot S$ is a $Q$-local martingale. Now $H\cdot S-Z^{Q}$ is a positive local 
$Q$-martingale and by Fatou lemma therefore also $Q$-supermartingale. Because
$Z^{Q}$ is a true martingale, $H\cdot S$ itself must be a supermartingale.
\end{proof}

In the next step we will construct a candidate optimal trading strategy and
prove that its wealth process is a supermartingale under any $Q\in \mathcal{M}\cap P_{\hat{V}}$. 
Note that the supermartingale property holds over the
larger set $P_{\hat{V}}$ rather than just those measures that lead to bliss-free expected utility $P_{V}$.

\begin{proposition}
\label{prop: strategy H}Assume there is $\bar{Q}\in \mathcal{M}\cap P_{\hat{V}}^{e}$. 
Under the assumptions of Proposition~\ref{prop: U(X_n) converges to
U(hatX) in L1(P)} there is a trading strategy $H\in L(S)$, a sequence of
maximizing tame strategies $H^{(n)}$ and a semimartingale $\tilde{V}$ such that
\begin{enumerate}
\item \label{P:H:item1}$\tilde{V}$ is $\bar{Q}$-supermartingale;

\item \label{P:H:item2}$\tilde{V}=\mathrm{rqlim}_{n\rightarrow \infty }H^{(n)}\cdot S,$ see
equation (\ref{eq: rqlim});

\item \label{P:H:item3}$H\cdot S\geq \tilde{V}$ and $H\cdot S-\tilde{V}$ is an increasing
process;

\item \label{P:H:item4}In particular, $H^{(n)}\cdot S_{T}\overset{P\text{-a.s.}}{\rightarrow }\tilde{V}_{T}$;

\item \label{P:H:item5}$U(B+H^{(n)}\cdot S_{T})\overset{L^{1}(P)}{\rightarrow }U(B+\tilde{V}_{T})$ 
and thus $I_{U}(B+\tilde{V}_{T})=u(B)\in \mathbb{R}$;

\item \label{P:H:item6}$H\cdot S$ is a $Q$-supermartingale for any $Q\in \mathcal{M}\cap P_{\hat{V}}$.
\end{enumerate}
\end{proposition}

\begin{proof}
\emph{Step 1)} First we prove that there is a maximizing sequence $\tilde{H}^{(n)}\in \scr{T}$ 
such that for any $Q\in \mathcal{S}\cap P_{\hat{V}}$ there is 
a $Q$-martingale $Z^{Q}$ with the property $\tilde{H}^{(n)}\cdot S\geq Z^{Q}$.
This is similar in spirit to \citet[Proposition 3.8]{biagini.cerny.11}
but there each $Q$ calls for a different subsequence whereas here the
maximizing (sub)sequence will be the same for all $Q$-s. Proposition \ref{prop: U(X_n) 
converges to U(hatX) in L1(P)} gives a maximizing sequence $\hat{H}^{(n)}\in \scr{T}$, 
a random variable $\hat{X}\in L^{0}(P)$, and a uniform lower bound $0\geq R\in L^1(P)$ 
such that $R\leq U(B+\hat{H}^{(n)}\cdot S_{T})\rightarrow U(B+\hat{X})$ in $L^{1}(P)$ and $P$-a.s. 

\emph{Step 2)} Let $\tilde{W}=U^{-1}(R)\leq 0$. By the Fenchel inequality,
\begin{equation*}
-E^{Q}[\lambda \tilde{W}]\leq E[\hat{V}(\lambda dQ/dP)]-E[U(\tilde{W})]=I_{\hat{V}}(\lambda dQ/dP)-E[R],
\end{equation*}
we conclude that $\tilde{W}$ is in $L^{1}(Q)$ for any $Q\in P_{\hat{V}}$.
The wealth process $\tilde{H}^{(n)}\cdot S$ is a $Q$-supermartingale for any 
$Q\in \mathcal{S}$ and hence 
\begin{equation}\label{eq: Htilde cdot S ge ZQ}
E_{t}^{Q}[B]+\tilde{H}^{(n)}\cdot S_{t}\geq E_{t}^{Q}[B+\tilde{H}^{(n)}\cdot
S_{T}]\geq E_{t}^{Q}[-(B+\tilde{H}^{(n)}\cdot S_{T})^{-}]\geq E_{t}^{Q}[\tilde{W}],  
\end{equation}
therefore $Z_{t}^{Q}=E_{t}^{Q}[\tilde{W}]-E_{t}^{Q}[B]$ yields the
lower bound announced in Step 1).

\emph{Step 3)} Now select a fixed $\bar{Q}$ in $\mathcal{M}\cap P_{\hat{V}}^{e}$ and 
apply Theorem~D in \cite{delbaen.schachermayer.98} to construct
processes $\tilde{V}$, $H$, and a sequence 
\begin{equation}\label{eq: tails Htilde}
H^{(n)}\in \mathrm{conv}(\tilde{H}^{(n)},\tilde{H}^{(n+1)},\ldots ),
\end{equation}
with the properties claimed in items~\eqref{P:H:item1}--\eqref{P:H:item4}. Item~\eqref{P:H:item5} follows from 
$\tilde{V}_{T}=\hat{X}$ and from the fact that $H^{(n)}$ is still a maximizing sequence.

\emph{Step 4)} The uniform lower bound $Z^{Q}$ for $\tilde{H}^{(n)}\cdot S$
obtained in (\ref{eq: Htilde cdot S ge ZQ}) also applies to $H^{(n)}\cdot S$
for $H^{(n)}$ from (\ref{eq: tails Htilde}) and therefore by item~\eqref{P:H:item2} also to 
$\tilde{V}$ because $Z^{Q}$ can be chosen right-continuous. Consequently $H\cdot S\geq Z^{Q}$ and 
by step 2)\ in the proof of Proposition~\ref{prop: tame strategies are super} $H\cdot S$ is a 
$Q$-supermartingale for every $Q\in \mathcal{M}\cap P_{\hat{V}}$.
\end{proof}

In the second step we will prove that the candidate strategy $H$ attains
the maximal utility $u(B)$. Therefore by item~\eqref{P:H:item6} of 
Proposition~\ref{prop: strategy H} the optimizer $H$ belongs to the supermartingale class of
strategies. It is readily seen that our approach simplifies and generalizes
the results of \cite{schachermayer.01,schachermayer.03}, in particular we
completely sidestep dynamic optimization arguments in the proof of the
supermartingale property, see also \cite{owen.zitkovic.09}. This is all the
more remarkable because the tools we use do not go beyond those pioneered by
Schachermayer and his co-authors in the run-up to \cite{schachermayer.03}.

\begin{table}[tbp] \centering

\begin{tabular}{lcccc}
&  & $M^{\hat{V}}=L^{\hat{V}}$ &  & $M^{\hat{V}}\subsetneq L^{\hat{V}}\smallskip $ \\ \hline\hline
\vspace{-0.9\baselineskip}
&  &  &  &  \\ 
{$M^{\hat{U}}=L^{\hat{U}}$} &  & none required; &  & $\exists\lambda >1;\lambda \hat{Y}\in \mathrm{dom}\,I_{V}$
\smallskip  \\ \hline
\vspace{-0.9\baselineskip}
&  &  &  &  \\ 
\multirow{3}{*}{$M^{\hat{U}}\subsetneq L^{\hat{U}}\nsim L^{\infty }$} 
&  & \multirow{3}{*}{$u(B)= \sup\limits_{X\in \mathcal{C}^{\circledast\circledast}}I_U(B+X)$;} 
&  & $u(B)= \sup\limits_{X\in \mathcal{C}^{\ast\ast}}I_U(B+X)$; \\ 
&  &  &  & \raisebox{0.15cm}{$\exists\lambda >1;\lambda \hat{Y}\in \mathrm{dom}\,I_{V}$}\smallskip\\ 
\hline \vspace{-0.9\baselineskip} &  &  &  &  \\ 
$M^{\hat{U}}\subsetneq L^{\hat{U}}\sim L^{\infty }$ &  & none required; &  & 
cannot occur$\smallskip $ \\ \hline\hline
\end{tabular}
\medskip 
\caption{Additional assumptions of Theorem~\ref{thm: main}.}\label{tab: assumptions}
\end{table}

\begin{theorem}
\label{thm: main}Assume

\begin{enumerate}
\item $B\in L^{\hat{U}}$;

\item $(B+\mathcal{C)}\cap \mathrm{core\,}\mathcal{D}\neq \emptyset $;

\item no arbitrage over $\mathcal{C}^{\circledast \circledast }$;

\item $\lim_{x\rightarrow \infty }\overline{u}(x)/x=0$;

\item $\mathcal{C}_{\sigma }^{\circledast }$ is norm-dense in $\mathcal{C}^{\circledast }$ 
(Assumption~\ref{ass: sig martingale measures are dense});

\item and further specific assumptions as detailed in Table \ref{tab:
assumptions}.
\end{enumerate}

Then Theorem~\ref{thm: strong duality}, Proposition~\ref{prop: U(X_n)
converges to U(hatX) in L1(P)}, Proposition~\ref{prop: y + Y} and
Proposition~\ref{prop: strategy H} apply and the strategy $H$ from
Proposition~\ref{prop: strategy H} is optimal and admissible.
\end{theorem}

\begin{proof}
\emph{Step 1)} No arbitrage over $\mathcal{C}^{\circledast \circledast }$
implies existence of an equivalent separating measure 
\citep[Proposition~3.5]{gao.xanthos.17}. Obtain 
$\bar{Q}\in \mathcal{M}\cap P_{\hat{V}}^{e}$ from 
Assumption~\ref{ass: sig martingale measures are dense}. We know from 
Proposition~\ref{prop: strategy H}, items~\eqref{P:H:item1}, \eqref{P:H:item3}, and \eqref{P:H:item5}, that 
$$I_{U}((B+H\cdot S_{T})\wedge 0)\geq
I_{U}((B+\tilde{V}_{T})\wedge 0)>-\infty .$$

\emph{Step 2)} By item \eqref{P:H:item6} of Proposition~\ref{prop: strategy H} one has 
$E[Y(H\cdot S_{T})]\leq 0$ for any $Y\in C_{\sigma }^{\circledast }$.
Therefore for any $m\ge0$ and any $Y\in C_{\sigma }^{\circledast }$ we obtain 
\begin{equation*}
\left\langle \left(B + H\cdot S_{T}\right) \wedge m,Y\right\rangle \leq
\left\langle B + H\cdot S_{T} ,Y\right\rangle \leq \langle B,Y\rangle.
\end{equation*}

\emph{Step 3)} By Proposition~\ref{prop: y + Y} and 
Assumption~\ref{ass: sig martingale measures are dense} there are sequences 
$W_{n}\in L^{\hat{V}}$ and $Y_{n}\in C_{\sigma }^{\circledast }$ with the 
properties $(Y_{n}+W_{n})\in \mathcal{C}^{\circledast }$, $\Vert W_{n}\Vert _{\hat{V}}\to 0$, and
$$I_{V}(Y_{n}+W_{n})+\left\langle B,Y_{n}+W_{n}\right\rangle \to u(B).$$
Step 1)\ implies $\Vert\!\left( B+H\cdot S_{T}\right)\wedge m\Vert _{\hat{U}}<\infty $ and the Fenchel inequality yields
\begin{align*}
I_{U}(\left(B+ H\cdot S_{T}\right) \wedge m) &\leq 
I_{V}(Y_{n}+W_{n})+\left\langle \left(B+ H\cdot S_{T}\right) \wedge
m,Y_{n}+W_{n}\right\rangle \\
&\leq I_{V}(Y_{n}+W_{n})+\langle B,Y_n+W_n\rangle\\
&\qquad\qquad+\langle \left(B+ H\cdot S_{T}\right) \wedge m - B,W_n\rangle,
\end{align*}
where we have used $\left\langle \left(B+ H\cdot S_{T}\right) \wedge m,Y_{n}\right\rangle \leq \langle B,Y_n\rangle$ from step 2).

\emph{Step 4)}  Use the H\"{o}lder inequality for Orlicz spaces and let $n\rightarrow \infty $ in step 3) 
to obtain $I_{U}(\left(B+ H\cdot S_{T}\right) \wedge m)\leq u(B)$. By monotone convergence \citep[Theorem~2.8.2]{bogachev.07}
letting $m\nearrow \infty $ we find $I_{U}(B+H\cdot S_{T})\leq u(B)$. 
Items~\eqref{P:H:item3} and \eqref{P:H:item5} of Proposition~\ref{prop: strategy H} now yield
\begin{equation}
I_{U}(B+H\cdot S_{T})=I_{U}(B+\tilde{V}_{T})=u(B).
\label{eq: I_U(H cdot S_T) = u}
\end{equation}

\emph{Step 5) }It remains to show that $H\cdot S$ can be approximated by $H^{(n)}\cdot S$ 
at intermediate times in the sense $\mathrm{rqlim}_{n\rightarrow \infty }H^{(n)}\cdot S= \tilde{V}=H\cdot S$ 
on the set $B+H\cdot S_{T}<\overline{x}$, recalling $\overline{x}$ in equation (\ref{eq: bliss point}). 
When $U$ is strictly monotone, the inequality $H\cdot S_{T}\geq \tilde{V}_{T}= \hat{X}$ 
together with equality (\ref{eq: I_U(H cdot S_T) = u}) imply 
\begin{equation}
H\cdot S_{T}=\tilde{V}_{T}.  \label{eq: H cdot S_T}
\end{equation}
By Proposition~\ref{prop: strategy H} the process $H\cdot S-\tilde{V}$ is
non-negative and increasing which in view of (\ref{eq: H cdot S_T}) is only
possible if $H\cdot S=\tilde{V}$.

\emph{Step 6)} When $\overline{x}<\infty $ argue by contradiction. Suppose
there is a non-null set $A$ on which $H\cdot S_{t}>\tilde{V}_{t}$ for some $t $ (not necessarily the same on each path) 
and $B+\tilde{V}_{T}<\overline{x}$. As $H\cdot S-\tilde{V}$ is increasing, it follows 
$H\cdot S_{T}>\tilde{V}_{T}$ on $A$ and because $B+\tilde{V}_{T}<\overline{x}$ on $A$ this
contradicts $I_{U}(B+\tilde{V}_{T})=I_{U}(B+H\cdot S_{T})$.
\end{proof}

\section{Conclusions}

\label{sect: conclusions}We have studied expected utility maximization from
the point of view of conjugate duality over Orlicz spaces $(L^{\hat{U}},L^{\hat{V}})$ 
determined by the left tail of the utility function and the right
tail of its conjugate, respectively. In this setup objects in $L^{\hat{V}}$
can be interpreted as complete market state price densities but not necessarily as
separating measures. In Theorem~\ref{thm: strong duality} we have
established the Fenchel duality over state price densities, applicable also to
large financial markets, in circumstances where none of the standard
regularity conditions apply. In Theorem~\ref{thm: main} we have provided
construction of the optimal trading strategy that does not rely on the dual
maximizer $\hat{Q}$ being a separating measure or being equivalent to $P$.
In the case $L^{\hat{U}}\sim L^{\infty }$ we have achieved this goal under
the minimal conditions from the seminal works of \citet{kramkov.schachermayer.03} 
and \citet{hugonnier.kramkov.04}.

The Fenchel duality formula mentioned above, 
\begin{equation*}
\sup_{X\in \mathcal{C}}I_{U}(X)=\min_{Y\in L^{\hat{V}}}\{I_{V}(Y)+\delta _{\mathcal{C}\cap \mathcal{D}}^{\circledast }(Y)\},
\end{equation*}
has an interesting economic interpretation. The quantity $\delta _{\mathcal{C}\cap \mathcal{D}}^{\circledast }(d\hat{Q}/dP)$ 
can be interpreted as an increase in the initial endowment required to bring the expected utility in a
complete market $\hat{Q}$ to the optimal level $u(0)$. In the same vein the
term $\delta _{\mathcal{C}\cap \mathcal{D}}^{\circledast }(dQ/dP)$ may be
interpreted as the shadow price of the implicit trading constraint presented
by the finiteness of the effective domain $\mathcal{D}$ of expected utility.
This provides a fresh perspective on the classical results of 
\citet{kramkov.schachermayer.99,kramkov.schachermayer.03}.

Our analysis motivates the study of the fundamental theorem of
asset pricing in an Orlicz space setting in a small financial market. The
question is whether $\mathcal{C}^{\star \star }\cap L_{+}^{\hat{U}}=\{0\}$
implies $\mathcal{C}=\mathcal{C}^{\circledast \circledast }$ and whether
the set of $\sigma $-martingale measures with density in $L^{\hat{V}}$ is
norm-dense among all separating measures, in the sense of Assumption~\ref{ass: sig martingale measures are dense}. 
The answer is known to be affirmative when $L^{\hat{U}}=L^{\infty }$ from the work of 
\citet{kabanov.97, delbaen.schachermayer.98}.

\pdfbookmark{Appendices}{appx} \appendix

\section{Key results in convex duality}

\label{sect: duality summary}In this appendix we have collected results from
convex analysis required in the main body of the paper, principally in the
proofs of Section~\ref{sect: duality over separating}. Unless explicitly
specified the functions are defined over a locally convex, Hausdorff
topological vector space $(E,\tau )$. Let $E^{\prime }$ denote the
topological dual of $(E,\tau )$, namely the space of linear, continuous
functionals on $(E,\tau )$. Any other topology $\sigma $ on $E$ such that
its topological dual coincides with $E^{\prime },$ meaning $E^{\prime
}=(E,\sigma )^{\prime }$ is called compatible with the dual pair 
\citep[Section~3]{rockafellar.74}. The conjugate functions are defined on 
$E^{\prime }$, endowed with a
topology $\tau ^{\prime }$ compatible with the dual pair, namely such that
the dual space of $(E^{\prime },\tau ^{\prime })$ equals $E$. Taking $E$ and 
$E^{\prime }$ as fixed it is known that the coarsest compatible topology on $E$ 
is the initial topology $\sigma (E,E^{\prime })$ while the finest
compatible topology in $E$ is the Mackey topology $\tau (E,E^{\prime })$.
For $y\in E^{\prime }$ and $x\in E$ we denote the bilinear form $y(x)$ by $\left\langle x,y\right\rangle $.

\begin{definition}
\label{def: effective domain}For a concave (resp. convex) function $h$ with
values in $[-\infty ,\infty ]$ its effective domain $\mathrm{dom\,}h$ is
defined by $\mathrm{dom\,}h=\{x:h(x)>-\infty \}$ (resp. $\{x:h(x)<\infty \}$). 
A concave (resp. convex) function $h$ is called proper if $h<\infty $
(resp. $h>-\infty $) and $\mathrm{dom\,}h$ is non-empty.
\end{definition}

\begin{definition}
\label{def: usc}A function $h$ (not necessarily concave/convex) with values
in $[-\infty ,\infty ]$ is called upper semi-continuous (resp. lower
semicontinuous), in short u.s.c. (resp. l.s.c.), if for each $c\in \mathbb{R}
$ the set $\{x:h(x)\geq c\}$ (resp. $\{x:h(x)\leq c\}$) is closed.
\end{definition}

\begin{definition}
\label{def: usc hull}For a function $h$ we denote by $\mathrm{usc}\mathsf{\,}h$ 
the upper semicontinuous hull of $h,$ i.e. the smallest upper
semicontinuous function that dominates $h$. Likewise, for a convex function $h$ we denote by $\mathrm{lsc}\mathsf{\,}h$ the lower semicontinuous hull of $h$, i.e. the greatest lower semicontinuous function dominated by $h$.
\end{definition}

\begin{proposition}
\label{prop: h = usc h}The upper semicontinuous hull (resp. the lower
semicontinuous hull) is given by the formula
\begin{equation*}
\mathrm{usc}\mathsf{\,}h(x)=\sup_{x_{\alpha }\rightarrow x}\underset{\alpha }{\lim \sup }\,h(x_{\alpha }),
\end{equation*}
resp. 
\begin{equation*}
\mathrm{lsc}\mathsf{\,}h(x)=\inf_{x_{\alpha }\rightarrow x}\underset{\alpha }{\lim \inf }\,h(x_{\alpha }),
\end{equation*}
where nets can be replaced by sequences when $(E,\tau )$ is first-countable
(in particular normed). A function $h$ is u.s.c. (resp. l.s.c.) if and only
if $h\geq \mathrm{usc}\mathsf{\,}h$ (resp. $h\leq \mathrm{lsc\,}h$).
\end{proposition}

\begin{proof}
\citet[eq.~3.7]{rockafellar.74} and \citet[Lemma~2.42]{aliprantis.border.06}.
\end{proof}

\begin{lemma}
\label{lem: sum of usc is usc}The sum of two u.s.c. functions with values in 
$[-\infty ,\infty )$ (resp. two l.s.c. functions with values in $(-\infty
,\infty ]$) is u.s.c. (resp. l.s.c.).
\end{lemma}

\begin{proof}
We have 
\begin{align*}
\mathrm{lsc}\,(g+h)(x) &=\inf_{x_{\alpha }\rightarrow x}\underset{\alpha }{\lim \inf }(g+h)(x_{\alpha }) \\
&\geq \inf_{x_{\alpha }\rightarrow x}\underset{\alpha }{\{\lim \inf }g(x_{\alpha })
+\underset{\alpha }{\lim \inf }h(x_{\alpha })\} \\
&\geq \inf_{x_{\alpha }\rightarrow x}\underset{\alpha }{\lim \inf }g(x_{\alpha })
+\inf_{x_{\alpha }\rightarrow x}\underset{\alpha }{\lim \inf }h(x_{\alpha }) \\
&=\mathrm{lsc}\mathsf{\,}g(x)+\mathrm{lsc}\mathsf{\,}h(x)=g(x)+h(x)
\end{align*}
and the statement follows by Proposition~\ref{prop: h = usc h}.
\end{proof}

\begin{theorem}
\label{thm: lsc hull invariance}For a concave (resp. convex) function $h$
one has 
\begin{equation*}
\mathrm{usc}_{\sigma (E,E^{\prime })}\mathsf{\,}h=\mathrm{usc}_{\tau
(E,E^{\prime })}\mathsf{\,}h
\end{equation*}
resp. $\mathrm{lsc}_{\sigma (E,E^{\prime })}\mathsf{\,}h=\mathrm{lsc}_{\tau
(E,E^{\prime })}\mathsf{\,}h$, meaning that the upper (resp. lower)
semicontinuous hull of a concave (resp. convex) function is the same in any
compatible topology.
\end{theorem}

\begin{proof}
See \cite{aliprantis.border.06}, Theorem~5.98 and Corollary~5.99.
\end{proof}

\begin{proposition}
\label{prop: usc h is proper}Suppose $h$ is concave (resp. convex). If $\mathrm{usc}\mathsf{\,}h$ 
(resp. $\mathrm{lsc}\mathsf{\,}h$) is finite-valued at a point then necessarily $\mathrm{usc}\mathsf{\,}h$ 
(resp. $\mathrm{lsc}\mathsf{\,}h$) is proper.
\end{proposition}

\begin{proof}
See \citet[Theorem~4]{rockafellar.74}.
\end{proof}

\begin{definition}[\protect\cite{rockafellar.74}]
\label{def: closure}For a concave function $h$ the upper closure $\cl h$ is defined 
\begin{equation*}
\cl h=\left\{ 
\begin{array}{cc}
\mathrm{usc}\mathsf{\,}h & \text{if }\mathrm{usc}\mathsf{\,}h<\infty \\ 
\infty & \text{otherwise}
\end{array}
\right..
\end{equation*}
Likewise, for a convex function $h$ the lower closure $\cl h$ is defined as
\begin{equation*}
\cl h=\left\{ 
\begin{array}{cc}
\mathrm{lsc}\mathsf{\,}h & \text{if }\mathrm{lsc}\mathsf{\,}h>-\infty \\ 
-\infty & \text{otherwise}
\end{array}
\right. \text{.}
\end{equation*}
We say that $h$ is closed if $h=\cl h$.
\end{definition}

\begin{proposition}
\label{prop: h*}For $h$ concave (resp. convex) $h^{\ast }=(\mathrm{usc}\mathsf{\,}h)^{\ast }$ 
(resp. $h^{\ast }=(\mathrm{lsc}\mathsf{\,}h)^{\ast }$) and $h^{\ast }=(\cl h)^{\ast }$ is closed.
\end{proposition}

\begin{proof}
See \citet[Theorem~5]{rockafellar.74} and \citet[Theorem~2.3.1]{zalinescu.02}.
\end{proof}

\begin{theorem}[Fenchel-Moreau]
\label{thm: Fenchel-Moreau}For $h$ concave or convex $h^{\ast \ast }=\cl h$.
\end{theorem}

\begin{proof}
See \citet[Theorem~5]{rockafellar.74}.
\end{proof}

\begin{corollary}
\label{cor: h* is proper}A concave (convex)\ function $\cl h$ is proper if and only if $h^{\ast }$ is proper.
\end{corollary}

\begin{proposition}
\label{prop: gao.xanthos.18}When $L^{\hat{U}}=(M^{\hat{V}})^{\star }$ a
concave (resp. convex) function $h$ on $L^{\hat{U}}$ is u.s.c. (resp.
l.s.c.) in the duality $(L^{\hat{U}},M^{\hat{V}})$ if and only if
\begin{align*}
h(x) &=\mathrm{usc}\mathsf{\,}h(x)=\sup_{\substack{ x_{n}\rightarrow xP\text{-a.s.}  \\ 
\left\Vert x_{n}\right\Vert _{\hat{U}}<K}}\underset{n}{\lim \sup }\,h(x_{n}),\text{ resp.} \\
h(x) &=\mathrm{lsc}\mathsf{\,}h(x)=\inf_{\substack{ x_{n}\rightarrow x  \\ 
\left\Vert x_{n}\right\Vert _{\hat{U}}<K}}\underset{n}{\lim \inf }\,h(x_{n}),
\end{align*}
for all $x\in L^{\hat{U}}$. That is, in computing a candidate for
u.s.c./l.s.c. hull in the $(L^{\hat{U}},M^{\hat{V}})$ duality nets can be
replaced with a.s.-convergent norm-bounded sequences.
\end{proposition}
\begin{proof}
\citet[Theorem~2.4]{gao.xanthos.18}.
\end{proof}

\begin{definition}
\label{def: inf convolution}Suppose $f,g$ are two concave and proper
functions. Their sup\-re\-mal convolution $f\conv g:E\rightarrow [
-\infty ,\infty ]$ is defined as 
\begin{equation*}
f\conv g(x)=\sup_{z\in E}\{f(x-z)+g(z)\}.
\end{equation*}
Likewise, for two proper convex functions $f,g$ their inf(imal) convolution
is given by 
\begin{equation*}
f\conv g(x)=\inf_{z\in E}\{f(x-z)+g(z)\}.
\end{equation*}
\end{definition}

\begin{lemma}
\label{lem: (f+g)*}For $f,g$ proper concave (convex) one has
\begin{equation}
\left( f\conv g\right) ^{\ast }=f^{\ast }+g^{\ast }.
\label{eq: (f square g)*}
\end{equation}
For concave $f$ and $g$ such that $f^{\ast }$ and $g^{\ast }$ are proper and 
$\cl f+\cl g=\cl (f+g)$ one has for all $y\in
E^{\prime }$
\begin{equation}
\cl\left( f^{\ast }\conv g^{\ast }\right)
(y)=(f+g)^{\ast }(y)=\inf_{x\in E}\{\left\langle x,y\right\rangle
-(f(x)+g(x))\}.  \label{eq: (f* square g*)}
\end{equation}
\end{lemma}

\begin{proof}
Formula (\ref{eq: (f square g)*}) follows from an easy computation 
\citep[eq. 9.30]{rockafellar.74}. The same formula applied to $f^{\ast }$ 
and $g^{\ast }$ yields
\begin{equation}
\left( f^{\ast }\conv g^{\ast }\right) ^{\ast }=f^{\ast \ast }+g^{\ast
\ast }=\cl f+\cl g,  \label{eq: (f* square g*)*}
\end{equation}
where the second equality follows by Theorem~\ref{thm: Fenchel-Moreau}. By
Proposition~\ref{prop: h*}
\begin{equation*}
\left( f+g\right) ^{\ast }=(\cl(f+g))^{\ast }
=\left( \cl f+\cl g\right) ^{\ast }
=\left( f^{\ast }\conv g^{\ast }\right)^{\ast \ast }
=\cl\left( f^{\ast }\conv g^{\ast }\right),
\end{equation*}
where the last two equalities follow from (\ref{eq: (f* square g*)*}) and
again Theorem~\ref{thm: Fenchel-Moreau}. The last equality in (\ref{eq: (f*
square g*)}) is immediate from the definition of conjugate function.
\end{proof}

\begin{proposition}
\label{prop: I_U is usc}When $U$ decreases superlinearly at $-\infty $ the
expected utility functional $I_{U}$ is $\ast $-upper semicontinuous. In the
linear case $I_{U}$ is $L^{1}$-norm continuous everywhere and therefore $\circledast $-u.s.c.
\end{proposition}

\begin{proof}
i) In the superlinear case $\lim_{x\rightarrow -\infty }U(x)/x=\infty $ one
has $(M^{\hat{V}})^{\star }=L^{\hat{U}}$. By Propositions~\ref{prop: h = usc
h} and \ref{prop: gao.xanthos.18} it suffices to prove that for every
pointwise convergent norm-bounded sequence $X_{n}\rightarrow X$ one has 
\begin{equation}
\underset{n\rightarrow \infty }{\lim \sup }\,I_{U}(X_{n})\leq I_{U}(X).  \label{eq: lim sup I_U(x_n)}
\end{equation}

Because $U(0)=0$ and $U$ is increasing and concave, Fatou lemma gives
\begin{align}
\underset{n\rightarrow \infty }{\lim \sup }\,I_{U}(-X_{n}^{-}) &\leq
I_{U}(-X^{-}),  \label{eq: usc estimate 1} \\
\underset{n\rightarrow \infty }{\lim \sup }\,\{I_{U}(X_{n}^{+})-U_{+}^{\prime }(0)E[X_{n}^{+}]\} 
&\leq I_{U}(X^{+})-U_{+}^{\prime }(0)E[X^{+}].\label{eq: usc estimate 2}
\end{align}
$M^{\hat{V}}$ with the Orlicz norm is an order-continuous Banach lattice. By \citet[Theorem~2.1]{gao.14}
$X_{n}$ is $\sigma (L^{\hat{U}},M^{\hat{V}})$-convergent to $X$. By \citet[Proposition~3.6]{wickstead.08}
the lattice operations are $\sigma (L^{\hat{U}},M^{\hat{V}})$-sequentially continuous on norm bounded subsets of $L^{\hat{U}}$ and hence $\lim_{n\rightarrow \infty }E[X_{n}^{+}]=E[X^{+}].$ On combining inequalities (\ref{eq: usc estimate 1}, \ref{eq: usc estimate 2}) we thus obtain (\ref{eq: lim sup I_U(x_n)}), which completes the proof.

ii) In the remaining linear case $\lim_{x\rightarrow -\infty }U(x)/x<\infty $
the space $L^{\hat{U}}$ is isomorphic to $L^{1}$ and $I_{U}$ is finite
everywhere. $I_{U}$ is norm-continuous on $L^{1}$ because $I_{U}$ is bounded below on any norm-bounded neighbourhood
of $0$, see \citep[Theorem~5.43]{aliprantis.border.06}. In this case the norm topology is
compatible with the duality and norm-continuity therefore implies $\circledast $-upper semicontinuity by 
Theorem~\ref{thm: lsc hull invariance}.
\end{proof}

\section{Corner solution with exponential utility}

\label{sect: appxB}

In this appendix we take $U(x)=-e^{-x}$. A routine calculation yields $V(y)=y\ln y-y$ and

\begin{equation}
\max_{y>0}I_{V}(yZ)=-E[Z]e^{-E[Z\ln Z]/E[Z]}.  \label{eq: max I_V(yZ)}
\end{equation}
The model for asset price $X$ and the optimal strategy are described in
Sections~\ref{sect: appxB asset price}-\ref{sect: appxB dual optimizer
sequence}. The non-existence of a supermartingale deflator with terminal
value $U^{\prime }(-X_{T})$, where $-X_{T}$ is the optimal terminal wealth
and $\hat{Y}=U^{\prime }(-X_{T})$ is the dual optimizer from 
Theorem~\ref{thm: strong duality}, is shown in Section~\ref{sect: no weak
supermartingale deflator}, where it is also noted that the optimal wealth
process $-X$ is a submartingale under $\hat{Q}$, $d\hat{Q}/dP=\hat{Y}/E[\hat{Y}]$.

\subsection{Asset price process}

\label{sect: appxB asset price}Let $X$ be a special semimartingale L\'{e}vy
process with characteristics $(b^{X},0,F^{X})$ where $b^{X}\in \mathbb{R}$
and
\begin{equation*}
F^{X}(dx)=\frac{3}{4\sqrt{\pi }}x^{-5/2}e^{-x}I_{(0,\infty )}(x)dx+\delta_{-1/2}(dx),
\end{equation*}
with $\delta _{x}$ denoting a Dirac measure at point $x$. Consequently the
cumulant generating function of $X$ is given by 
\begin{align*}
\kappa _{X}(v) &=b^{X}v+\int (e^{vx}-1-vx)F^{X}(dx) \\
&=e^{-v/2}+\left( 1-v\right) ^{3/2}-2+(2+b^{X})v\text{ for }v\leq 1\text{
and }\infty \text{ otherwise.}
\end{align*}
$X$ can be interpreted as a sum of a compensated one-sided (positive)
tempered stable process with parameters $\beta =3/2$, $\alpha =1/\Gamma
(-\beta )$, $\lambda =1$; a compensated Poisson process with intensity $1$
and jump size $-1/2$; and a drift component with drift $b^{X}$.

In this construction it is important that $\beta >1$. The choice of the tempered
stable process for positive jumps is significant only to the extent that its
L\'{e}vy measure density is exponential divided by a polynomial of
sufficiently high degree as $x\rightarrow \infty $; any other L\'{e}vy
measure with this property would do just as well. The convenience of
the tempered stable formulation is that it yields a simple expression for the
cumulant generating function \citep{kuechler.tappe.13}
which makes it particularly obvious that we will be dealing with a corner
solution.

The choice of the Poisson process for the single negative jump is not important,
but the jump size being bounded below by $-1/2$ means that $\mathcal{E}(X)$
is strictly positive and so our example could be recast in terms of an
exponential L\'{e}vy model. We will not pursue this line of exposition here
and instead formulate everything as trading on $X$.

To this end, $\kappa _{X}(v)$ being finite for $v\leq 1$ and exponential
being a submultiplicative function \citep[Proposition~25.4]{sato.99},
we obtain $\sup_{t\in [ 0,T]}\left\vert X_{t}\right\vert \in L^{\hat{U}}$ \citep[Theorem~28.18]{sato.99}. 
By \citet[Proposition~6.4]{biagini.cerny.11}
this means every separating measure in $L^{\hat{V}}$ is a local martingale
measure for $X$ and by Proposition~\ref{prop: tame strategies are super}
every separating measure is therefore a supermartingale measure for all $L^{\hat{U}}$-tame strategies.

\subsection{Candidate optimal trading strategy}

Consider now optimization over buy-and-hold strategies in $X$. Assume $X_{0}=0$ so that 
the terminal wealth reads $\vartheta X_{T}$. Expected utility
is then $I_{U}(\vartheta X_{T})=-\exp \left( \kappa _{X}(-\vartheta \right))$. 
Optimization over $\vartheta $ yields the following first order
condition for an interior maximum, 
\begin{equation*}
0=\kappa _{X}^{\prime }(-v)=-\frac{1}{2}e^{v/2}-\frac{3}{2}\left( 1+v\right)
^{1/2}+2+b^{X}.
\end{equation*}
Provided $b^{X}<-2+1/(2\sqrt{e})\approx -1.7$, which is what we assume
hereafter, there will be no interior optimizer and instead maximum will be
achieved at $\vartheta =-1$. For future reference let
\begin{equation*}
\kappa _{X}^{\prime }(1)=b^{X}+\int x\left( e^{x}-1\right)
F^{X}(dx)=b^{X}+2-1/(2\sqrt{e})=:-A<0.
\end{equation*}
Our task is to prove that $-X_{T}$ is the optimal wealth and therefore $\vartheta = -1$ is 
the optimal strategy and $-\exp \left( \kappa^{X}(1\right) )$ is the maximal utility. 
The buy-and-hold strategy $\vartheta =-1$ is $L^{\hat{U}}$-tame in view of 
$\sup_{t\in [ 0,T]}\left\vert X_{t}\right\vert \in L^{\hat{U}}$.

\subsection{Dual optimizing sequence of separating measures}

\label{sect: appxB dual optimizer sequence}For $n=1,2,\ldots $ define 
\begin{align*}
W_{n}(x) &=e^{x}-1+\frac{4\sqrt{\pi }}{3}x^{5/2}e^{x}K_{n}1_{[n,n+1]}(x) \\
&= e^{x}-1+K_{n}1_{[n,n+1]}(x)dx/F^{X}(dx), \\
K_{n} &=A/(n+1/2).
\end{align*}
With this definition one has 
\begin{equation}
b^{X}+\int xW_{n}(x)F^{X}(dx)=\underbrace{b^{X}+\int \left( e^{x}-1\right)
xF^{X}(dx)}_{-A}+(n+1/2)K_{n}=0.  \label{eq: zero drift}
\end{equation}
Let $J^{X}$ be the jump measure associated with the process $X$ and define
\begin{equation}
d\mathcal{L}(Z^{(n)})=dX+(W_{n}(x)-x)dJ^{X}.  \label{eq: dL(Z)}
\end{equation}

Note that $Z^{(n)}$ is a well-defined strictly positive process because
jumps on the right-hand side are bounded below by $1/\sqrt{e}-1$. It follows
that the $P$-drift of $\mathcal{L}(Z^{(n)})$ is given by
\begin{equation}
b^{\mathcal{L}(Z^{(n)})}=b^{X}+\int (W_{n}(x)-x)F^{X}(dx),
\label{eq: b_L(Z)}
\end{equation}
and that $Z^{(n)}\exp (-b^{\mathcal{L}(Z^{(n)})}t)$ is a density process of
a local martingale measure for $X$ by virtue of the Girsanov theorem and (\ref{eq: zero drift}) which may be rewritten as
\begin{equation*}
b^{X}dt+d\langle
X,\mathcal{L}(Z^{(n)}\exp (-b^{\mathcal{L}(Z^{(n)})}t))\rangle=0.
\end{equation*}
Denote this local martingale measure by $Q^{(n)}$, with 
$$dQ^{(n)}=Z_{T}^{(n)}\exp (-b^{\mathcal{L}(Z^{(n)})}T)dP.$$ 
Here $\left\langle
X,Y\right\rangle $ now stands for the predictable quadratic covariation, i.e.
the drift part of process $[X,Y]$ (provided the semimartingale $[X,Y]$ is
special). By the It\={o} formula
\begin{align}
b^{\ln Z^{(n)}} &=b^{X}+\int \left( \ln (1+W_{n}(x))-x\right) F^{X}(dx) 
\notag \\
&=b^{X}+\underbrace{\int_{n}^{n+1}\ln (1+K_{n}\frac{dx}{F(dx)}e^{-x})F(dx)}
_{B_{n}\searrow 0}.  \label{eq: b_lnZ}
\end{align}

From the Girsanov theorem the drift of $\ln Z^{(n)}$ under $Q^{(n)}$ is
given by
\begin{align*}
b_{Q^{(n)}}^{\ln Z^{(n)}} &=b^{\ln Z^{(n)}}+\int W_{n}(x)\ln
(1+W_{n}(x))F^{X}(dx), \\
&=b^{\ln Z^{(n)}}+\int W_{n}(x)\ln \left(e^{x}+K_{n}\frac{dx}{F^{X}(dx)}
1_{[n,n+1]}(x)\right)F^{X}(dx) \\
&=b^{X}+B_{n}+\int W_{n}(x)xF^{X}(dx)\\
&\qquad+\underbrace{\int_{n}^{n+1}W_{n}(x)\ln
\left(1+K_{n}\frac{dx}{F^{X}(dx)}e^{-x}\right)F^{X}(dx)}_{C_{n}\searrow 0} \\
&=B_{n}+C_{n}\searrow 0,
\end{align*}
where we have substituted for $b^{\ln Z^{(n)}}$ from (\ref{eq: b_lnZ}) and
used (\ref{eq: zero drift}) in the penultimate line.

Apply the complete market utility formula (\ref{eq: u_Q Fenchel}), 
$u_{Q^{(n)}}(0)=\max_{y>0}I_{V}(yZ_{T}^{(n)})$, and evaluate it using the
expression (\ref{eq: max I_V(yZ)}), and the help of identities 
$E[Z_{T}^{(n)}]=\exp (b^{\mathcal{L}(Z^{(n)})}T)$ and $E[Z_{T}^{(n)}\ln
Z_{T}^{(n)}]/E[Z_{T}^{(n)}]=E^{Q^{(n)}}[\ln Z_{T}^{(n)}]=b_{Q^{(n)}}^{\ln
Z^{(n)}}T$, 
\begin{equation*}
u_{Q^{(n)}}(0)=-E[Z_{T}^{(n)}]e^{-E[Z_{T}^{(n)}\ln
Z_{T}^{(n)}]/E[Z_{T}^{(n)}]}=-\exp (b^{\mathcal{L}(Z^{(n)})}T-b_{Q^{(n)}}
^{\ln Z^{(n)}}T).
\end{equation*}
In (\ref{eq: b_L(Z)}) substitute for $W_{n}$ and rearrange to obtain
\begin{equation*}
b^{\mathcal{L}(Z^{(n)})}=b^{X}+\int \left( e^{x}-1-x\right)
F^{X}(dx)+\int_{n}^{n+1}K_{n}dx\searrow \kappa _{X}(1).
\end{equation*}
In conclusion, 
\begin{equation*}
u_{Q^{(n)}}(0)=-\exp ((b^{\mathcal{L}(Z^{(n)})}-b_{Q^{(n)}}^{\ln
Z^{(n)}})T)\searrow -\exp (\kappa ^{X}(1))=E[U(-X_{T})],
\end{equation*}
which proves optimality of $\vartheta =-1$ as all tame strategies are
supermartingales under the measures $Q^{(n)}$ and hence the utility of any
tame strategy may not exceed the expression on the left-hand side.

\subsection{There is no supermartingale deflator ending with 
\texorpdfstring{$U^{\prime}(-X_{T})$}{U'(X\textunderscore T)}}

\label{sect: no weak supermartingale deflator}The distribution $P_{X_{t}}$
is absolutely continuous with respect to the Lebesgue measure 
\citep[Lemma~27.1, Theorem~27.7]{sato.99}, and its support is the entire real 
line. We wish to investigate whether
there is $c\geq 0$ and supermartingale $D$ with 
$$D_{T}=U^{\prime}(-X_{T})=e^{X_{T}}$$ 
such that $D(c-X)$ is also a supermartingale. Now $X$
is a process with independent increments so one can evaluate the conditional
expectations $E_{t}[(c-X_{T})e^{X_{T}}]$ explicitly.

Take $W(x)=e^{x}-1$ then $Z_{T}$ defined by (\ref{eq: dL(Z)}) gives
precisely $Z_{T}= e^{X_{T}}$ and we can reuse the calculations in the
previous section with $K_{n}= 0$ and $d\hat{Q}/dP=e^{X_{T}-\kappa
_{X}(1)T}=\hat{Y}/E[\hat{Y}]$ to obtain
\begin{align*}
E_{t}[(c-X_{T})e^{X_{T}}] &=e^{X_{t}+\kappa ^{X}(1)(T-t)}\left(
c-X_{t}+E_{t}[(X_{t}-X_{T})e^{X_{T}-X_{t}-\kappa ^{X}(1)(T-t)}]\right) \\
&=e^{X_{t}+\kappa ^{X}(1)(T-t)}\left( c-X_{t}-b_{\hat{Q}}^{X}(T-t)\right)
\leq D_{t}(c-X_{t}).
\end{align*}
Now $-b_{\hat{Q}}^{X}= A>0$ and so when $c<X_{t}<c+A,$ which happens
with non-zero probability thanks to the support of $X$ being the whole real
line, the left-hand side is positive while the right-hand side, no matter
how one chooses $D_{t}\geq 0$, is non-positive. Therefore, $U^{\prime
}(-X_{T})$ cannot be identified with a terminal value of even a weak
supermartingale deflator.

 The inequality 
$-b_{\hat{Q}}^{X}= A>0$ implies that $-X$ is a $\hat{Q}$-submartingale because $X$ is a L\'{e}vy process under $\hat{Q}$.

\section{Utility may increase from \texorpdfstring{$\mathcal{C}$}{C} to 
\texorpdfstring{$\mathcal{C}^{\circledast\circledast}$}{C**}}

\label{sect: utility increases C to C**}

We have observed in Theorem~\ref{thm: strong duality} that the maximal
utility over $\mathcal{C}$ and its norm-closure $\mathcal{C}^{\star \star }$
always coincide. Thus when the norm topology is compatible with the economic
duality we are guaranteed that the maximal utility over $\mathcal{C}$ and 
$\mathcal{C}^{\circledast \circledast }$ is the same. We study equivalent
conditions for such compatibility in Theorem~\ref{thm: topology equivalences}
.

In Section~\ref{sect: utility increases} we then provide an explicit 
\emph{arbitrage-free} example, using exponential utility, of a situation where
maximal utility increases by going from $\mathcal{C}$ to 
$\mathcal{C}^{\circledast \circledast }$. This construction can be applied to
arbitrary utility in the setting of statement 8) in Theorem~\ref{thm:
topology equivalences}, that is on any probability space that is not purely
atomic and such that the Orlicz heart $M^{\hat{U}}$ is strictly contained in
the Orlicz space $L^{\hat{U}}.$ We thus find that the difficult cases are
essentially those where $M^{\hat{U}}\neq L^{\hat{U}}$. The word essentially
refers to infinite, purely atomic spaces where it is not known whether in
all instances with $M^{\hat{U}}\neq L^{\hat{U}}$ such an example exists, the
same way it is not known whether 7) implies 8) in Theorem~\ref{thm: topology
equivalences}.

\subsection{When is the norm topology compatible with economic duality?}

\label{sect: topology equivalences}

It turns out that the characterization hinges on the properties of a
functional called `modular'. We therefore begin with a more general concept
of an ordered modular space, following the exposition of \citet{nowak.89a}. 
We then specialize this more general setup to Orlicz spaces used in this
paper. The attribute $\sigma $ should be read as \textquotedblleft
countably\textquotedblright\ or \textquotedblleft
countable\textquotedblright .

Let $E$ be a $\sigma $-Dedekind complete Riesz space. A functional 
$\rho :E\rightarrow [ 0,\infty ]$ is called a \emph{modular} if the following
conditions hold:

(i) $\rho (x)=0$ if and only if $x=0$.

(ii) $\left\vert x\right\vert <\left\vert y\right\vert $ implies $\rho
(x)<\rho (y)$.

(iii) $\rho (x_{1}\vee x_{2})<\rho (x_{1})+\rho (x_{2})$ for $x_{1}\geq
0,x_{2}\geq 0$.

(iv) $\rho (\lambda x)\rightarrow 0$ if $\lambda \rightarrow 0$.

\noindent One can verify that with this definition $\rho $ is a modular also
in the original sense of \citet{musielak.orlicz.59}. A modular $\rho $ 
is said to be convex, if $\rho (\lambda_{1}x_{1}
+\lambda _{2}x_{2})<\lambda _{1}\rho (x_{1})+\lambda _{2}\rho(x_{2})$ 
for $\lambda _{1},\lambda _{2}\geq 0$ and $\lambda _{1}+\lambda
_{2}=1$. A modular $\rho $ is said to be \emph{metrizing} whenever $\rho
(x_{n})\rightarrow 0$ implies $\rho (2x_{n})\rightarrow 0$ for a sequence 
$\{x_{n}\}$ in $E$. Recall the definition of the corresponding gauge norm, 
$\left\Vert x\right\Vert _{\rho }=\inf \{\lambda >0:\rho (x/\lambda )\leq 1\}$.  
For the definition of modular topology see \citet[p. 262]{nowak.89a}.

The following characterization of modular and norm convergence is key. A
sequence $\{x_{n}\}$ in $E$ converges to zero modularly if and only if there is 
$\lambda >0$ such that $\rho (\lambda x_{n})\rightarrow 0,$ while it
converges to zero in the\ norm $\left\Vert \cdot \right\Vert _{\rho }$ if and only if 
$\rho (\lambda x_{n})\rightarrow 0$ for all $\lambda >0$, \emph{ibid}.
Therefore $\rho $ fails to be metrizing precisely when there is a sequence
that converges to zero modularly but not in the norm.

Let $\Phi $ be a Young function, $(\Omega ,\mathcal{F},P)$ a probability
measure space. Then 
$\rho :L^{0}(\Omega ,\mathcal{F},P)\rightarrow [ 0,\infty ]$
\begin{equation}
\rho (X)=E[\Phi (X)],  \label{eq: Orlicz rho}
\end{equation}
is a convex orthogonally additive modular on the Orlicz space $L^{\Phi }$,
satisfying the $\sigma $-Lebesgue property, the $\sigma $-Fatou property and
the $\sigma $-Levi property (\citealp[Section~2]{nowak.89a}; \citealp[pp.~274-275]{nowak.89b}). 
The norm $\left\Vert \cdot \right\Vert _{\rho }$ is known as the
Luxemburg norm in this setting.

Finally, recall the definition of Mackey topology and strong topology for a
given dual pair \citep[Sections 5.18 and 5.19]{aliprantis.border.06}. The following 
theorem gives full characterization of circumstances under
which the norm closure and the economic closure coincide.

\begin{theorem}
\label{thm: topology equivalences}Let $\Phi $ be a Young function, $\rho $
the corresponding modular from (\ref{eq: Orlicz rho}) and $\Psi $ the
conjugate of $\Phi .$ Then the strong topology $\beta (L^{\Phi },L^{\Psi })$
coincides with the norm topology on $L^{\Phi },$ the Mackey topology $\tau
(L^{\Phi },L^{\Psi })$ coincides with the modular topology and the following
are equivalent:

\begin{enumerate}
\item $\rho $ is metrizing;

\item Every sequence $\{X_{n}\}\in L^{\Phi }$ modularly convergent to zero
is also norm-convergent;

\item $\beta (L^{\Phi },L^{\Psi })=\tau (L^{\Phi },L^{\Psi })$;

\item The gauge norm $\left\Vert \cdot \right\Vert _{\rho }$ is
order-continuous on $L^{\Phi }$;

\item $(L^{\Phi },\tau (L^{\Phi },L^{\Psi }))$ is barreled;

\item $(L^{\Phi })^{\star }=L^{\Psi },$ that is $L^{\Psi }$ is the norm-dual
of $L^{\Phi }$;

If, furthermore, we exclude the case where the probability space is finite
(i.e. we assume $P$ is \emph{not }supported on a finite number of atoms)
then the following is equivalent to 1)-6)

\item $M^{\Phi }=L^{\Phi }$.

If we also exclude the case where $P$ is purely atomic then the following is
equivalent to 1)-7)

\item $\Phi $ satisfies so-called $\Delta _{2}$--condition at $\infty$, i.e. 
there is $K>0$ and $x_{0}>0$ such that $\Phi (2x)\leq K\Phi (x)$ for all $x>x_{0}$.
\end{enumerate}

The implications 8) $\Rightarrow $ 7) $\Rightarrow $ 1)-6) hold without
further assumptions.
\end{theorem}

\begin{proof}
\citet[Theorems~3.2 and 4.2]{nowak.89b}
show that $\beta (L^{\Phi },L^{\Psi })$ is the norm topology on $L^{\Phi }$
and $\tau (L^{\Phi },L^{\Psi })$ is the modular topology. The equivalence 
$1)\iff 2)$ is trivial. Equivalences $1)\iff 3)\iff 4)$ follow from 
\citet[Theorem~2.3]{nowak.89a}, while $3)\iff 5)$ follows from a standard 
result in topology, \citet[Corollary~II.1.2 and II.2.4]{husain.khaleelulla.78}. 
Equivalence of $4)\ $and $6)$ for $\Phi $ finite follows from 
\citet[Ch~15, p.~336 and Ch~19, pp.~572-3]{zaanen.83}, while in the remaining 
case $L^{\Phi }\sim L^{\infty },L^{\Psi }\sim L^{1}$ both 4)\ and 6) are true 
if the probability is finite and both are false otherwise.

$7)\Rightarrow 6)$ follows from \citet[Theorem~2.2.11]{edgar.sucheston.92}
for finite $\Phi $ while for $\Phi $ that jumps to infinity $7)$ is false
and the implication holds trivially.

$2)\Rightarrow 7)$ distinguish two cases: A) When $L^{\Phi }\sim L^{\infty }$
use non-finiteness of the probability space to construct a disjoint sequence
of events $A_{n}\in \mathcal{F},P(A_{n})>0$ and $\sum_{n=1}^{\infty
}P(A_{n})=1$. Without loss of generality we may assume $\Phi (1)<\infty
,\Phi (2)=\infty .$ Let $B_{n}=\bigcup\limits_{k=1}^{n}A_{k}$ and define 
$X_{n}(\omega )=1_{B_{n}^{c}}(\omega ).$ Then $I_{\Phi }(X_{n})\rightarrow 0$
by dominated convergence while $I_{\Phi }(2X_{n})=\infty $ meaning 
$\{X_{n}\} $ converges to zero modularly but not in norm. This shows
statement $2)$ is false and the implication holds trivially. B) When $\Phi $
is everywhere finite argue by contradiction. Suppose there is $X\in L^{\Phi
}\setminus M^{\Phi }$ which in particular means $X\notin L^{\infty }$ and
therefore $(\Omega ,\mathcal{F},P)$ must be a non-finite probability space.
Without loss of generality we may suppose $0\leq X,$ $I_{\Phi }(X)<\infty
,I_{\Phi }(2X)=\infty $. Define $X_{n}=X1_{X\geq n}$. Once again $I_{\Phi
}(X_{n})\rightarrow 0$ by dominated convergence while $I_{\Phi
}(2X_{n})=\infty $ in contradiction to $2).$

The implication $8)\Rightarrow 7)$ follows from \citet[Theorem~131.3]{zaanen.83}. 
The opposite implication for $P$ that is not purely atomic follows from \citet[Theorem~III.2]{rao.ren.91}.
\end{proof}

\subsection{Illustrative example}

\label{sect: utility increases} Let $U(x)=-e^{-x}$ which implies $V(y)=y\ln
y-y$. Consider a probability space $(\Omega =\mathbb{Z},\mathcal{F}=2^{
\mathbb{Z}},P)$ with $P(\{n\})=e^{-\left\vert n\right\vert }$ for 
$n\in \{1,\pm 2,\pm 3,\ldots \},$ $P(\{-1\})=e^{-5},$ and 
\begin{equation*}
P(\{0\})=1-\sum_{\left\vert n\right\vert \geq 1}P(\{n\})=1-2(e^{2}-
e)^{-1}-e^{-1}-e^{-5}.
\end{equation*}

Define a random variable $X$ by setting $X(-1)=-1,X(1)=1,$ and $X=0$
elsewhere. Let $Y(n)=n$ for $\left\vert n\right\vert \geq 2$ and $Y=0$
otherwise. Define a sequence of random variables $\{X_{k},Y_{k}\}_{k\in 
\mathbb{N}}$ by setting $Y_{k}=Y1_{\left\vert Y\right\vert \geq k}$ and 
$X_{k}=X+Y_{k}$. Note that $Y,$ in common with all $Y_{k}$, has finite
exponential moments in the interval $(-1,1)$ but not beyond. This means 
$\{Y_{k}\}$ converges to zero modularly (which by previous theorem means in the
economic duality) but not in the norm on $L^{\hat{U}}$. By construction
\begin{equation}
E[Y_{k}|X]=0\text{ for all }k\in \mathbb{N}\text{.}  \label{eq: E[Y|X]}
\end{equation}

Think of $X_{k}$ as an excess return on a traded position. Define the
marketed subspace as 
$\mathcal{K}=\mathrm{span}(\{X_{k}\}_{k\in \mathbb{N}})$. The probability 
measure $\hat{Q}$ defined by $$d\hat{Q}/dP=e^{-2X}/E[e^{-2X}]$$
 is a bliss-free completion of the market as,
by construction, $E^{\hat{Q}}[X]=E^{\hat{Q}}[Y_{k}]=0$ for $k\in \mathbb{N}$
and $\hat{Q}$ has finite entropy, 
\begin{equation*}
H(\hat{Q}||P)=E\left[ \frac{d\hat{Q}}{dP}\ln \frac{d\hat{Q}}{dP}\right]
=-\ln E[e^{-2X}].
\end{equation*}
One can readily verify that $2X$ is the optimal wealth in the complete
market $\hat{Q}$, because for exponential utility and any bliss-free state
price measure $Q$ one has by formula (\ref{eq: u_Q Fenchel}) and direct
calculation
\begin{equation*}
u_{Q}(0)=\sup_{E^{Q}[X]\leq
0}I_{U}(X)=\min_{y>0}I_{V}(ydQ/dP)=I_{V}(e^{-H(Q||P)}dQ/dP)=-e^{-H(Q||P)}.
\end{equation*}

However, the maximal utility in the original market $\mathcal{K}$ is
strictly lower than $u_{Q}$. To see this, consider a finite linear
combination $Z\in \mathcal{K}$ with $1\leq k(1)<k(2)<\ldots <k(N)$ being the
indices in ascending order of vectors with non-zero coefficients,
\begin{align*}
Z &=\sum_{i=1}^{N}\lambda _{i}X_{k(i)}=\xi _{N}X_{k(N)}+\sum_{j=1}^{N-1}\xi
_{N-j}\left( X_{k(N-j)}-X_{k(N-j+1)}\right) \\
&=\xi _{N}(X+Y_{k(N)})+\sum_{j=1}^{N-1}\xi _{N-j}\left(
Y_{k(N-j)}-Y_{k(N-j+1)}\right) ,
\end{align*}
where $\xi _{N-j}=\left( \sum_{i=1}^{N-j}\lambda _{i}\right) $ for 
$j=0,\ldots ,N-1$. The random variables $X$ and $\{Y_{k(N-j)}-Y_{k(N-j+1)}\}_{j=1}^{N-1}$ 
are in $L^{\infty }$. It follows 
$$E[e^{-Z}]<\infty \iff E[\exp \left( -\xi _{N}Y_{k(N)}\right) ]\iff E[\exp \left( -\xi
_{N}Y\right) ]\iff \left\vert \xi _{n}\right\vert <1.$$ 
In view of (\ref{eq: E[Y|X]}) the conditional Jensen's inequality 
\citep[Lemma~2.1]{mussmann.88} yields 
\begin{equation*}
E[e^{-Z}]\geq E[e^{-\xi _{N}X}]\geq E[e^{-X}],
\end{equation*}
where the last inequality holds because
$E[e^{-\lambda X}]=e^{-5}e^{\lambda }+e^{-1}e^{-\lambda }$ is a strictly
convex function of $\lambda $ attaining its global minimum at $\lambda =2$ and
therefore decreasing on $(-\infty ,2].$ It follows that that the maximal
utility over $\mathcal{K}$ satisfies 
\begin{equation*}
u(0)=\sup_{W\in \mathcal{K}}I_{U}(W)=\sup_{W\in \mathcal{C}}I_{U}(W)=-E[e^{-X}].
\end{equation*}

Finally, let us examine the economic closure $\mathcal{C}^{\circledast
\circledast }$. As $\{Y_{k}\}$ converges to $0$ in the economic duality, we
have $\lambda X\in \mathcal{C}^{\circledast \circledast }$ for all 
$\lambda \in \mathbb{R}$. In contrast, 
\begin{equation*}
\lambda X\in \mathrm{cl}^{\circledast }(\mathcal{C}\cap \mathcal{D})
\text{ if and only if }\left\vert \lambda \right\vert \leq 1.
\end{equation*}
Consequently,
\begin{align*}
-E[e^{-X}] &= \sup_{W\in \mathcal{C}}I_{U}(W)\\
&<\max_{W\in \mathcal{C}^{\circledast \circledast }}I_{U}(W)
=\min_{Y\in \mathcal{C}^{\circledast}}I_{V}(Y)=-e^{-H(\hat{Q}||P)}=-E[e^{-2X}].
\end{align*}
Note, however, that Theorem~\ref{thm: strong duality} continues to hold for 
$\mathcal{A}=\mathcal{C}$ 
(as well as for $\mathcal{A}=\mathcal{C}^{\circledast \circledast }$ 
which is just the last three equalities above):
\begin{align*}
-E[e^{-X}] &= \sup_{W\in \mathcal{C}}I_{U}(W)=\max_{W\in \mathrm{cl}
^{\circledast }(\mathcal{C}\cap \mathcal{D})}I_{U}(W)=\min
I_{V}(Y)+\sup_{W\in \mathcal{C}\cap \mathcal{D}}E[WY] \\
&=I_{V}(e^{-X})+\sup_{W\in \mathcal{C}\cap \mathcal{D}}
E[We^{-X}]=I_{V}(e^{-X})+E[Xe^{-X}],
\end{align*}
and the optimal effective completion of $\mathcal{C}$ is given by the state
price density $$d\hat{Q}/dP=e^{-X}/E[e^{-X}].$$ Note that under $\hat{Q}$ the
optimal wealth process \emph{increases} in expectation, $E^{\hat{Q}}[X]>0$, hence 
$X$ is a $\hat{Q}$--submartingale.

In the given example $X\in L^{\infty }\subseteq M^{\hat{U}}$ possesses all
exponential moments and the dual optimizer is therefore a separating
measure. One can modify this example by splitting the state $\{0\}$ into
countably many states where $X$ is unbounded and such that it only possesses
exponential moment of order at most, say, 1.5 while maintaining the present
inequality $E[Xe^{-1.5X}]>0$. In this way one may exhibit a situation where 
\begin{equation*}
\max_{W\in \mathrm{cl}^{\circledast }(\mathcal{C}\cap \mathcal{D})}
I_{U}(W)<\max_{W\in \mathcal{C}^{\circledast \circledast }}I_{U}(W),
\end{equation*}
the first optimizer is $X,$ the second optimizer is $1.5X$ and each
optimizer represents a corner solution. 

In the first case the corner is
caused by the `nuisance' zero-mean shocks $Y_{k}$ which do not allow us to
increase our position in $X$ beyond $1$ while we are trading inside $\mathcal{C}$. 
These nuisance shocks `stop contaminating' $X$ as one passes
to the economic closure $\mathcal{C}^{\circledast \circledast }$. One is
now able to take a position $\lambda X$ with $\lambda $ above $1$. In the
second case the corner over $\mathcal{C}^{\circledast \circledast }$ at 
$\lambda =1.5$ is inherent in $X$ itself. Seen in this light, duality over
separating measures (\ref{eq: duality gap}) signifies that even this corner
can be `removed' by passing to a full completion whose utility is
arbitrarily close to that of $\mathcal{C}^{\circledast\circledast }$.

\end{document}